\documentclass[MSNbibl,number,citesort,dvips]{arxbj}
\usepackage{upgreek}
\usepackage{graphicx}

\aid{0}
\volume{20}
\issue{2}
\pubyear{2014}
\firstpage{514}
\lastpage{544}
\doi{10.3150/12-BEJ496} 

\makeatletter
\newcommand{\RMo}{\mathrm{o}}
\newcommand{\RMO}{\mathrm{O}}

\newcommand{\cal}{\mathcal}

\newcommand{\R}{\mathbb{R}}
\newcommand{\E}{\mathbb{E}}
\renewcommand{\L}{\mathbb{L}}
\renewcommand{\P}{\mathbb{P}}
\renewcommand{\hat}{\widehat}
\renewcommand{\d}{\mathrm{d}}

\newtheorem{theo}{Theorem}[section]
\newtheorem{prop}{Proposition}[section]
\newtheorem{corol}{Corollary}[section]
\newtheorem{lemm}{Lemma}[section]

\makeatother

\begin{document}
\begin{frontmatter}

\title{Testing monotonicity via local least concave~majorants}
\runtitle{Testing monotonicity via local least concave majorants}

\begin{aug}
\author[1]{\fnms{Nathalie} \snm{Akakpo}\thanksref{1}\ead[label=e1]{nathalie.akakpo@upmc.fr}},
\author[2]{\fnms{Fadoua} \snm{Balabdaoui}\thanksref{2}\ead[label=e2]{fadoua@ceremade.dauphine.fr}} \and
\author[3]{\fnms{C\'{e}cile} \snm{Durot}\corref{}\thanksref{3}\ead[label=e3]{cecile.durot@gmail.com}}
\runauthor{N. Akakpo, F. Balabdaoui and C. Durot} 
\address[1]{Univ. Pierre et Marie Curie, CNRS UMR 7599, LPMA, 4 place
Jussieu, 75252 Paris Cedex 05, France. \printead{e1}}
\address[2]{Univ. Paris-Dauphine, Ceremade, Place du Mar\'echal De
Lattre De Tassigny, 75775 Paris Cedex 16, France. \printead{e2}}
\address[3]{Univ. Paris-Ouest, UFR SEGMI, 200 avenue de la r\'
{e}publique, 92001 Nanterre Cedex, France.\\ \printead{e3}}
\end{aug}

\received{\smonth{7} \syear{2011}}
\revised{\smonth{6} \syear{2012}}

%
\begin{abstract}
We propose a new testing procedure for detecting localized departures
from monotonicity of a signal embedded in white noise. In fact, we
perform simultaneously several tests that aim at detecting departures
from concavity for the integrated signal over various intervals of
different sizes and localizations. Each of these local tests relies on
estimating the distance between the restriction of the integrated
signal to some interval and its least concave majorant. Our test can be
easily implemented and is proved to achieve the optimal uniform
separation rate simultaneously for a wide range of H\"olderian
alternatives. Moreover, we show how this test can be extended to a
Gaussian regression framework with unknown variance. A simulation study
confirms the good performance of our procedure in practice.
\end{abstract}

%

\end{frontmatter}

\section{Introduction}
Suppose that we observe on the interval $[0,1]$ a stochastic process
$F_n$ that is governed by the white noise model
%
\begin{equation}
\label{eqmodel} F_n(t)=\int_0^tf(x)
\,\d x+\frac{\sigma}{\sqrt n}W(t),
\end{equation}
where $n \ge1$ is a given integer, $\sigma>0$ is known, $f\dvtx [0,1]\to
\R$
is an unknown function on $[0,1]$ assumed to be integrable, and $W$ is
a standard Brownian motion on $[0,1]$ starting at 0. We aim to test the
null hypothesis that $f$ is non-increasing on $[0,1]$
against the general alternative that it is not.

Several non-parametric procedures have already been proposed, either in
model (\ref{eqmodel}) or, most often, in a regression model. To avoid
going back and forth between the two directions of monotonicity, we
will only talk about non-increasing monotonicity since the other
direction can be treated similarly via a straightforward
transformation. Without being exhaustive, we now review some tests that
are well suited in the regression framework for detecting global
departures from monotonicity. Bowman, Jones and Gijbels \cite{bowman98}
propose a procedure based on the smallest bandwidth for which a
kernel-type estimator of the regression function is monotone. Their
test rejects monotonicity if this critical bandwidth is too large.
Durot \cite{durot03} exploits the equivalence between monotonicity of a
continuous regression curve $g$ defined on $[0,1]$ and concavity of
$G\dvtx t\mapsto\int_0^t g(x) \,\d x, t \in[0,1]$. In the uniform design
setting, Durot's test rejects monotonicity when the supremum distance
between an empirical estimator of $G$ and its least concave majorant is
too large. The test has the correct asymptotic level, and has the
advantage of being easy to implement. Dom{\'\i}nguez-Menchero, Gonz{\'
a}lez-Rodr{\'\i}guez and L{\'o}pez-Palomo \cite{dominguez05} propose a
test statistic based on the $\L_2$ distance of a regression estimator
to the set of all monotone curves. Tailored for fixed design regression
models, their method allows for repeated measurements at a given design
point, and presents also the merit of being easy to implement. Still in
the fixed design setting, Baraud, Huet and Laurent \cite{bhl03}
construct a multiple test that, roughly speaking, rejects monotonicity
if there is at least one partition of $[0,1]$ into intervals, among a
given collection, such that the estimated projection of the regression
function on the set of piecewise constant functions on this partition
is too far from the set of monotone functions. They show that their
test has the correct level for any given sample size and study its
uniform separation rate, with respect to an $\L_2$ criterion, over H\"
olderian balls of functions. We may also mention some procedures for
testing monotonicity in other frameworks. For instance, Durot \cite
{durot08} or Groeneboom and Jongbloed \cite{groe11A,groe11B} test the
monotonicity of a hazard rate, whereas Delgado and Escanciano \cite
{delgado} test the monotonicity of a conditional distribution. As
Durot \cite{durot03}, they consider a statistic based on some distance
between an empirical estimator of the function of interest and its
least concave majorant.

Other procedures were considered to detect local departures from
monotonicity. Hall and Heckman \cite{hallheckman00} test negativity of
the derivative of the regression function via a statistic based on the
slopes of the fitted least-squares regression lines over small blocks
of observations. Gijbels, Hall, Jones and Koch \cite{gijbelskoch00}
propose two statistics based on signs of differences of the response
variable. In the case of a uniform fixed design and i.i.d. errors with
a bounded density, the authors study the asymptotic power of their
procedure against local alternatives under which the regression
function departs linearly from the null hypothesis at a certain rate.
Ghosal, Sen and van der Vaart \cite{ghosal00} test negativity of the
first derivative of the regression function via a locally weighted
version of Kendall's tau. The asymptotic level of their test is
guaranteed and their test is powerful provided that the derivative of
the regression function locally departs from the null at a rate that is
fast enough. D\"umbgen and Spokoiny \cite{ds01} and Baraud, Huet and
Laurent \cite{bhl05} both propose two multiple testing procedures,
either in the Gaussian white noise model for \cite{ds01} or in a
regression model for \cite{bhl05}. The former authors consider two
procedures based on the supremum, over all bandwidths, of kernel-based
estimators for some distance from $f$ to the null hypothesis. The
distance they estimate is either the supremum distance from $f$ to the
set of non-increasing functions or the supremum of $f'$. The latter
authors propose a procedure based on the difference of local means and
another one based on local slopes, a method that is akin to that of
Hall and Heckman \cite{hallheckman00}. In both papers, the uniform
separation rate of each test is studied over a range of H\"olderian
balls of functions. Last, let us mention the (unfortunately
non-conservative) alternative approach developed by Hall and Van
Keilegom \cite{hallkeilegom} for testing local monotonicity of a
hazard rate.

In this paper, we propose a multiple testing procedure that may be seen
as a localized version of the test considered by \cite{durot03}. Based
on the observation of $F_n$ in (\ref{eqmodel}), it rejects
monotonicity if there is at least one subinterval of $[0,1]$, among a
given collection, such that the local least concave majorant of $F_n$
on this interval is too far from $F_n$. Our test has the correct level
for any given $n$, independently of the degree of smoothness of~$f$.
Its implementation is easy and does not require bootstrap nor any
a priori choice of smoothing parameter. Moreover, we show that
it is powerful simultaneously against most of the alternatives
considered in \cite{ghosal00}. We also study the uniform separation
rate of our test, with respect to two different criteria, over a range
of H\"olderian balls of functions. We recover the uniform separation
rates obtained in \cite{ds01,bhl05}, as well as new ones, and check
that our test achieves the optimal uniform separation rate
simultaneously over the considered range of H\"olderian balls. Besides,
we are concerned with the more realistic Gaussian regression framework
with fixed design and unknown variance. We describe how our test can be
adapted to such a framework, and prove that it still enjoys similar
properties. Finally, we briefly discuss how our method could be
extended to more general models.

The organization of the paper is as follows. In Section \ref
{seclevel}, we describe our testing procedure and the critical region
in the white noise model for a prescribed level. Section \ref
{secpower} contains theoretical results about the power of the test.
In Section \ref{secregression}, we turn to the regression framework.
Section~\ref{secsimulation} is devoted to the practical implementation
of our test -- both in the white noise and regression frameworks -- and
to its comparison with other procedures via a simulation study. In
Section \ref{secdiscussion}, we discuss possible extensions of our
test to more general models. All proofs are postponed to Section \ref
{secallproofs} or to the supplementary material \cite{CFNsupp}.

\section{Testing procedure}\label{seclevel}
Let us fix a closed sub-interval $I$ of $[0,1]$ and denote by ${\cal
D}^I$ the set of all functions from $[0,1]$ to $\R$ which are
non-increasing on $I$. We start with describing a procedure for testing
the null hypothesis
\[
H_0^I\dvtx  f\in{\cal D}^I
\]
against $
H_1^I\dvtx  f\notin{\cal D}^I$
within the framework (\ref{eqmodel}). Since the cumulative function
%
\begin{equation}
\label{eqcumul} F(t)=\int_0^tf(x) \,\d x,\qquad t
\in[0,1],
\end{equation}
is concave on $I$ under $H_0^I$ and $F_n$ estimates $F$, we reject
$H_0^I$ when $F_n$ is ``too far from being concave on $I$''.
Thus, we consider a procedure based on a local least concave majorant.
For every
continuous function $G\dvtx [0,1]\to\R$, we denote by $\hat G^I$ the least
concave majorant
of the restriction of $G$ to $I$, and simply denote $\hat G^{[0,1]}$ by
$\hat G$. Our test statistic is then
\[
S_{n}^I=\sqrt{\frac{n}{\sigma^2|I|}}\sup_{t\in I}
\bigl(\hat F_n^I(t)-F_n(t) \bigr),
\]
where $|I|$ denotes the length of $I$, and we reject $H_0^I$ when
$S_n^I$ is too large. From the following lemma, $S_n^I$ is the supremum
distance between $\hat F^I_n$ and $F_n$ over $I$, normalized in such a
way that its distribution does not depend on $I$ under the least
favorable hypothesis that $f$ is constant on $I$.
%
\begin{lemm}\label{lemleastfav}
Let $I$ be a closed sub-interval of $[0,1]$. If $f\equiv c$ over $I$
for some $c\in\R$, then
$S_n^I=\sup_{t\in I} (\hat W^I(t)-W(t) )/\sqrt{|I|}$
and is distributed as
%
\begin{equation}
\label{eqZ} Z:=\sup_{t\in[0,1]} \bigl(\hat W(t)-W(t) \bigr).
\end{equation}
\end{lemm}
For a fixed $\alpha\in(0,1)$, we reject the null hypothesis $H_0^I$ at
level $\alpha$ if
%
\begin{equation}
\label{eqIregion} S_n^I>q(\alpha),
\end{equation}
where $q(\alpha)$ is calibrated under the hypothesis that $f$ is
constant on $I$, that is, $q(\alpha)$ is the $(1-\alpha)$-quantile of
$Z$. As stated in the following theorem, this test is of non-asymptotic
level $\alpha$ and the hypothesis that the function $f$ is constant on
$I$ is least favorable, under no prior assumption on~$f$.
%
\begin{theo}\label{theoIlevel}
For every $\alpha\in(0,1)$,
\[
\sup_{f\in{\cal D}^I}\P_f \bigl[S_n^I>q(
\alpha) \bigr]=\alpha.
\]
Moreover, the supremum is achieved at $f\equiv c$ over $I$, for any
fixed $c\in\R$.
\end{theo}

We now turn to our main goal of testing the null hypothesis
%
\begin{equation}
\label{eqH0}
H_0\dvtx  f\in{\cal D}
\end{equation}
against the general alternative $ H_1\dvtx  f\notin{\cal D}$, where ${\cal
D}$ denotes the set of non-increasing functions from $[0,1]$ to $\R$.
In the case where the monotonicity assumption is rejected, we would
also like to detect the places where the monotonicity constraint is not
satisfied. Therefore, we consider a finite collection ${\cal C}_n$ of
sub-intervals of $[0,1]$, that may depend on $n$, and propose to
combine all the local tests of $H_0^I$ against $H_1^I$ for $I\in{\cal
C}_n$. In the spirit of the heuristic union-intersection principle of
Roy \cite{Roy} (see Chapter 2), we accept $H_0$ when we accept $H_0^I$
for all $I\in{\cal C}_n$, which leads to\vadjust{\goodbreak} $\max_{I\in{\cal C}_n} S_n^I$
as a natural test statistic. More precisely, given $\alpha\in(0,1)$, we
reject the null hypothesis $H_0$ at level $\alpha$ if
%
\begin{equation}
\label{eqregion}
\max_{I\in{\cal C}_n} S_n^I>s_{\alpha,n},
\end{equation}
where $s_{\alpha,n}$ is calibrated under the hypothesis that $f$ is a
constant function, that is, $s_{\alpha,n}$ is the $(1-\alpha
)$-quantile of
\[
\max_{I\in{\cal C}_n}\sqrt{\frac{1}{|I|}}\sup_{t\in I} \bigl(\hat
W^I(t)-W(t) \bigr).
\]
If $H_0$ is rejected, then we are able to identify one or several
intervals $I\in{\cal C}_n$ where the monotonicity assumption is
violated: these are intervals where $S_n^I>s_{\alpha,n}$.
Moreover, Theorem \ref{theolevel} below shows that this multiple
testing procedure has a non-asymptotic level~$\alpha$.
%
\begin{theo}\label{theolevel}
For every $\alpha\in(0,1)$,
\[
\sup_{f\in\cal D}\P_f \Bigl[\max_{I\in{\cal C}_n}
S_n^I>s_{\alpha,n} \Bigr]=\alpha.
\]
Moreover, the above supremum is achieved at $f\equiv c$ on $[0,1]$ for
any fixed $c\in\R$.
\end{theo}
%
Let us recall that the distribution of $S_n^I$ does not depend on $I$
under the least favorable hypothesis. This allows us to perform all the
local tests with the same critical threshold $s_{\alpha,n}$, which is
of practical interest. On the contrary, the multiple test in \cite
{bhl05} involves several critical values, which induces a complication
in its practical implementation (see~Section 5.1 in \cite{bhl05}).

The least concave majorants involved in our procedure can be computed
by first
discretizing the intervals $I\in{\cal C}_n$ and then using, for
example, the Pool Adjacent Violators Algorithm on a finite number of
points, see \cite{BBBB}, Chapter 1, page 13. Thus, $s_{\alpha,n}$ can be
computed using Monte Carlo simulations and the test is easily
implementable. However, in the white noise model this requires
simulating a large number of (discretized) Brownian motion paths and
the computation of
a local least concave majorant on each $I\in{\cal C}_n$, which can be
computationally
expensive. An alternative to Monte Carlo simulations is to replace the
critical threshold $s_{\alpha,n}$ with an upper-bound that is easier to
compute. Two different proposals for an upper-bound are given in the
following lemma.
%
\begin{lemm}\label{lemmajq}
For every $\gamma\in(0,1)$, let $q(\gamma)$ denote the $(1-\gamma
)$-quantile of the variable $Z$ defined in (\ref{eqZ}). Then, for
every $\alpha\in(0,1)$, we have
\[
s_{\alpha,n}\leq q \biggl(\frac{\alpha}{|{\cal C}_n|} \biggr)\leq
2\sqrt{2\log
\biggl(\frac{2|{\cal C}_n|}{\alpha} \biggr)}.
\]
\end{lemm}
Consider the test that rejects $H_0$ if
%
\begin{equation}
\label{eqregion2}
\max_{I\in{\cal C}_n} S_n^I>t_{\alpha,n},\vadjust{\goodbreak}
\end{equation}
where $t_{\alpha,n}$ denotes either $q (\alpha/|{\cal
C}_n| )$
or $2\sqrt{2\log(2|{\cal C}_n|/\alpha)}$. Combining Theorem
\ref{theolevel} with Lem\-ma~\ref{lemmajq} proves that this test is at
most of non-asymptotic level $\alpha$. In practice, the first proposal
of definition of $t_{\alpha,n}$ can be computed by using either recent
results of Balabdaoui and Pitman \cite{balabdpitman09} about a
characterization of the distribution function of $Z$ or fast Monte
Carlo simulations. Moreover, it is easy to see that all the theoretical
results we obtain in Section \ref{secpower} about the performance of
the test with critical region (\ref{eqregion}) continue to hold for
the test with critical region (\ref{eqregion2}) with both proposals
for $t_{\alpha,n}$. In practice however, the test with critical
region (\ref{eqregion2}) may have lower power than the test (\ref
{eqregion}), see Section \ref{secsimulation}.

\section{Performance of the test}\label{secpower}

Keeping in mind the connection between the white noise model and the
regression model of Section \ref{secregression}, we study the
theoretical performance of our test in model (\ref{eqmodel}) for
%
\begin{equation}
\label{eqCn} {\cal C}_n= \biggl\{ \biggl[\frac{i}{n},
\frac{j}{n} \biggr], i<j \mbox{ in }\{0,\ldots,n\} \biggr\},
\end{equation}
which can be viewed as the collection of all possible subintervals of
$[0,1]$ in model (\ref{eqregmod}). Of course, if we knew in advance the
interval $I$ over which $f$ is likely to violate the non-increasing
assumption, then the power would be largest for the choice ${\cal
C}_n=\{I\}$, which would be the right collection for testing $H_{0}^I$
instead of $H_{0}$. However, such a situation is far from being
realistic since $I$ is in general unknown. With the choice (\ref
{eqCn}), we expect ${\cal C}_n$ to contain an interval close enough
to~$I$. Therefore, by performing local tests simultaneously on all
intervals of ${\cal C}_n$, the resulting multiple testing procedure is
expected to detect a wide class of alternatives. Other possible choices
of ${\cal C}_{n}$ will be discussed in Section \ref{secdiscussion}.

In the sequel, we take $n\geq2$, and $\alpha$ and $\beta$ are some
fixed numbers in $(0,1)$. We will give a sufficient condition for our
test to achieve a prescribed power: we provide a condition on $f$ which
ensures that
%
\begin{equation}
\label{eqpower} \P_f \Bigl[\max_{I\in{\cal C}_n} S_n^I>s_{\alpha,n}
\Bigr]\geq1-\beta.
\end{equation}
Then, based on this condition, we study the uniform separation rate of
our test over H\"olderian classes of functions.



\subsection{Power of the test}\label{subsecCSpower}
For every $x<y$ in $[0,1]$, let
%
\begin{equation}
\label{eqbarf} \bar f_{xy}=\frac{1}{y-x}\int_x^y
f(u) \,\d u.
\end{equation}
%
In the following theorem, we prove that our
test achieves a prescribed power provided there are $x<t<y$ such that
$\bar f_{xy}$ is too large as compared to $\bar f_{xt}$.
%
\begin{theo}\label{theopowerslope}
There exists $C(\alpha,\beta)>0$ only depending on $\alpha$ and
$\beta$
such that (\ref{eqpower}) holds provided there exist $x,y\in[0,1]$
such that $y-x\geq2/n$ and
%
\begin{equation}
\label{eqpowerslope} \sup_{t\in[x,y]}\frac{t-x}{\sqrt{y-x}}(\bar f_{xy}-
\bar f_{xt})\geq C(\alpha,\beta)\sqrt\frac{\sigma^2\log n}{n}.
\end{equation}
\end{theo}

To illustrate this theorem, let us consider a sequence of alternatives where
$f$, which may depend on $n$, is assumed to be continuously
differentiable. Moreover, we assume that there exist an interval
$[t_n-\Delta_n,t_n+\Delta_n]\subset[0,1]$ and positive numbers $M$,
$\lambda_n$ and $\delta_n\leq\Delta_n$ such that
%
\begin{equation}
\label{eqpowerghosal} \cases{ f'(t)\geq0, &\quad on $[t_n-
\Delta_n,t_n+\Delta_n]$,
\cr
f'(t)\geq M\lambda_n, &\quad on $[t_n-
\delta_n,t_n+\delta_n]$.}
\end{equation}
%
Here, $t_n$, $\lambda_n$, $\delta_n$ and $\Delta_n$ may depend on $n$
while $M$ does not. Under these assumptions, we obtain the following corollary.
%
\begin{corol}\label{corolpowerslope}
There exists a positive number $C(\alpha,\beta)$ only depending on
$\alpha$ and $\beta$ such that (\ref{eqpower}) holds provided
$\Delta_n\geq1/n$ and
%
\begin{equation}
\label{eqcorolpowerslope} \frac{M\delta_n^2\lambda_n}{\sqrt{\Delta
_n}}\geq C(\alpha,\beta)\sqrt
\frac{\sigma^2\log n}{n}.
\end{equation}
\end{corol}
Corollary \ref{corolpowerslope} allows us to compare our test to the
one constructed in \cite{ghosal00} for monotonicity of a mean function
$f$ in a regression model with random design and $n$ observations.
Indeed, Theorem 5.2 in \cite{ghosal00}, when translated to the case of
non-increasing monotonicity as in (\ref{eqH0}) provides sufficient
conditions on $(\lambda_n,\delta_n,\Delta_n)$ for the test in \cite
{ghosal00} to be powerful against an alternative of the form (\ref
{eqpowerghosal}) with a large enough $M$. One can check that if
$(\lambda_n,\delta_n,\Delta_n)$ satisfies one of the sufficient
conditions given in Theorem 5.2 in \cite{ghosal00}, then it also satisfies
%
\begin{equation}
\label{eqcsghosal} \frac{\delta_n^2\lambda_n}{\sqrt{\Delta_n}}\geq
C\sqrt\frac{\sigma^2\log(1/\Delta_n)}{n}
\end{equation}
for some constant $C>0$ which does not depend on $n$. Hence, it also
satisfies (\ref{eqcorolpowerslope}) provided that
%
\begin{equation}
\label{eqcomparghosal} C\sqrt{\log(1/\Delta_{n})}\geq
\frac{C(\alpha,\beta)}{M}\sqrt{\log n}.
\end{equation}
It follows that Corollary \ref{corolpowerslope} shows that our test is
powerful simultaneously against all the alternatives considered in
Theorem 5.2 in \cite{ghosal00} for which $\Delta_n\geq1/n$ and
(\ref{eqcomparghosal}) holds. Noticing that (\ref{eqcomparghosal})
holds,
for instance, when $\Delta_n\leq n^{-\epsilon}$ for some $\epsilon>0$
and $M$ large enough, we conclude that our test is powerful
\textit{simultaneously} for most of the alternatives against which the test
in \cite{ghosal00} is shown to be powerful. However, as opposed to our
test the test in \cite{ghosal00} is \textit{not} simultaneously powerful
against those alternatives: a given alternative is detected only if a
parameter $h_n$ is chosen in manner that depends on that particular alternative.

As a second illustration of Theorem \ref{theopowerslope}, consider the
alternative that $f$ is U-shaped, that is $f$ is convex and there
exists $x_{0}\in(0,1)$ such that $f$ is non-increasing on $[0,x_{0}]$
and increasing on $[x_{0},1]$. Thus, $f$ violates the null hypothesis
on $[x_{0},1]$ and deviation from non-increasing monotonicity is
related to $f(1)-f(x_{0})>0$.
%
\begin{corol}\label{corolpowerUshaped}
Assume $f$ is U-shaped, denote by $x_{0}$ the greatest location of the
minimum of $f$, and let $\rho=f(1)-f(x_{0})$ and
\[
R=\lim_{x\uparrow1}\frac{f(x)-f(1)}{x-1}<\infty.
\]
Then, there exist a constant $C_{0} > 0$ and a real number
$C(\alpha,\beta) > 0$ only depending on $\alpha$ and $\beta$ such that
(\ref{eqpower}) holds provided $\rho>C_{0}R/n$ and
\[
\rho>C(\alpha,\beta) \biggl(\frac{\sigma^2R\log n}{n} \biggr)^{1/3}.
\]
\end{corol}

\subsection{Uniform separation rates}\label{subsecunifrates}
In this section, we focus on H\"olderian alternatives: we study the
performance of our test assuming that for some $R>0$ and $s\in(0,2]$
both unknown, $f$ belongs to the H\"olderian class ${\cal F}(s,R)$
defined as follows. For all $R>0$ and $s\in(0,1]$, ${\cal F}(s,R)$ is
the set of functions $f\dvtx [0,1]\to\R$ satisfying
\[
\bigl|f(u)-f(v)\bigr|\leq R|u-v|^{s}\qquad\mbox{for all }u,v \in[0,1],
\]
while for all $R>0$ and $s\in(1,2]$, ${\cal F}(s,R)$ is defined as the
set of differentiable functions $f\dvtx [0,1]\to\R$ with derivative $f'$ satisfying
$f'\in{\cal F}(s-1,R)$.

To evaluate the performance of our test, we consider a criterion that
measures the discrepancy of $f$ from $\cal D$. In the case where we
only assume that $f\in{\cal F}(s,R)$ for some $R>0$ and $s\in(0,2]$, we
consider the criterion
\[
\Delta_1(f)=\inf_{g\in\cal D}\sup_{t\in[0,1]}\bigl|f(t)-g(t)\bigr|,
\]
which is the supremum distance from $f$ to $\cal D$. Then, we restrict
our attention to the case where $f\in{\cal F}(s,R)$ for some $R>0$ and
$s\in(1,2]$
and consider the criterion
\[
\Delta_2(f)=\sup_{t\in[0,1]}f'(t).
\]
We recall that, for a given class of functions $\cal F$ and a given
criterion $\Delta$, the uniform separation rate of an $\alpha$-level
test $\Phi$ over $\cal F$ with respect to $\Delta$ is defined by
\[
\rho(\Phi,{\cal F},\Delta)=\inf\bigl\{\rho>0, \P_f(\Phi\mbox{
rejects }H_0)\geq1-\beta\mbox{ for all }f\in{\cal F}\mbox{ s.t. }
\Delta(f)\geq\rho\bigr\}
\]
and allows us to compare the performance of $\alpha$-level tests
: the smaller the better. The following theorem
provides an upper-bound for the uniform separation rate of our test
with respect to the criteria introduced above.
%
\begin{theo}\label{theoallpower}
Let ${\cal C}_n$ be the collection (\ref{eqCn}), $T_n$ be the test
with critical region (\ref{eqregion}), and let $R>0$ and $s\in(0,2]$.
Assume $n^s\sqrt{\log n}\geq R/\sigma$ and, in the case where $s\in
(1,2]$, assume moreover that $R\geq2^{1+2s}\sigma\sqrt{(\log n)/n}$.
Then, there exists a positive real $C(s,\alpha,\beta)$ only depending
on $s$, $\alpha$ and $\beta$ such that
%
\begin{equation}
\label{eqpowerbhl} \rho\bigl(T_n,{\cal F}(s,R),
\Delta_1\bigr)\leq C(s,\alpha,\beta)R^{1/(1+2s)} \biggl(
\frac{\sigma^2\log n}{n} \biggr)^{s/(1+2s)}
\end{equation}
and, in case $s\in(1,2]$,
%
\begin{equation}
\label{eqpowerds} \rho\bigl(T_n,{\cal F}(s,R),\Delta_2
\bigr)\leq C(s,\alpha,\beta)R^{3/(1+2s)} \biggl(\frac{\sigma^2\log
n}{n}
\biggr)^{(s-1)/(1+2s)}.
\end{equation}
\end{theo}
It should be noticed that the conditions
\[
n^s\sqrt{\log n}\geq R/\sigma\quad\mbox{and}\quad R\geq2^{1+2s}\sigma
\sqrt{(\log n)/n}
\]
simply mean that $n$ is sufficiently large when compared to $R/\sigma$.

Let us now discuss the optimality of the upper-bounds given in
Theorem \ref{theoallpower}. Indeed, Proposition \ref{proplowerbound}
below proves that, for each criterion $\Delta_1$ and $\Delta_2$ and any
choice of the smoothness parameter $s$ considered in Theorem \ref
{theoallpower}, no test can achieve a better rate over $\mathcal
{F}(s,R)$, up to a multiplicative constant that does not depend on $n$.
%
\begin{prop}\label{proplowerbound}
Assume $f\in\mathcal{F}(s,R)$ for some $R>0$ and $s\in(0,2]$. Let
$\mathcal{T}_\alpha$ be the set of all $\alpha$-level tests for testing
that $f\in\mathcal{D}$. Assume moreover that $R\geq\sigma\sqrt{(\log
n)/n^{1-\varepsilon}}$ for some $\varepsilon\in(0,1)$ and $n$ is
large enough. Then, there exists $\kappa(s,\alpha,\beta,\varepsilon)>0$
only depending on $s$, $\alpha$, $\beta$ and $\varepsilon$ such that
%
\begin{equation}
\label{eqlowerbound1} \inf_{\Phi_n \in\mathcal{T}_\alpha}\rho\bigl(\Phi
_n,{\cal
F}(s,R),\Delta_1\bigr)\geq\kappa(s,\alpha,\beta,
\varepsilon)R^{1/(1+2s)} \biggl(\frac
{\sigma^2\log n}{n} \biggr)^{s/(1+2s)}
\end{equation}
and, in case $s\in(1,2]$,
%
\begin{equation}
\label{eqlowerbound2} \inf_{\Phi_n\in\mathcal{T}_\alpha}\rho\bigl(\Phi
_n,{\cal
F}(s,R),\Delta_2\bigr)\geq\kappa(s,\alpha,\beta,\varepsilon)R^{3/(1+2s)}
\biggl(\frac
{\sigma^2\log n}{n} \biggr)^{(s-1)/(1+2s)}.
\end{equation}
\end{prop}
The lower-bound (\ref{eqlowerbound1}) is proved in \cite{bhl05}
(Proposition 2), and the other one in \cite{CFNsupp}. According
to (\ref{eqpowerbhl}) and (\ref{eqlowerbound1}), our test thus
achieves the optimal uniform separation rate (up to a multiplicative
constant) with respect to $\Delta_1$, simultaneously over all classes
${\cal F}(s,R)$ for $s\in(0,2]$ and a wide range of values of $R$. A
testing procedure enjoying similar performance in a Gaussian regression
model, at least for $s\in(0,1]$, can be found in \cite{bhl05}, Section~2.
In the white noise model we consider here, \cite{ds01} (Section 3.1)
propose a procedure that achieves the precise optimal uniform
separation rate (we mean, with the optimal constant) with respect to
$\Delta_1$, simultaneously over all classes ${\cal F}(1,R)$ with $R>0$.
But, to our knowledge, our test is the first one to achieve the
rate (\ref{eqpowerbhl}) in the case $s\in(1,2]$. On the other hand,
according to (\ref{eqpowerds}) and (\ref{eqlowerbound2}), our test
also achieves the optimal uniform separation rate (up to a
multiplicative constant) with respect to $\Delta_2$, simultaneously
over all classes ${\cal F}(s,R)$ for $s\in(1,2]$ and a wide range of
values of $R$. The second procedure proposed by \cite{bhl05} (Section
3) achieves this rate in a Gaussian regression model, and the second
procedure proposed by \cite{ds01} (Section 3.2) in the white noise
model is proved to achieve the optimal rate for $f\in{\cal F}(2,R)$.

\section{Testing procedure in a regression framework}\label{secregression}
In this section, we explain how our testing procedure can be extended,
with similar performance, to the more realistic regression model
%
\begin{equation}
\label{eqregmod} Y_i = f (i/n ) + \sigma\epsilon_i,\qquad
i=1,\ldots,n,
\end{equation}
where $f\dvtx [0,1]\to\R$ and $\sigma>0$ are unknown, and $
(\epsilon_i )_{1\leq i\leq n}$ are independent standard Gaussian variables.
Based on the observations $(Y_i)_{1\leq i\leq n}$, we would like to
test $H_0\dvtx  f\in{\cal D}$ against $H_1\dvtx f\notin{\cal D}$ where as above,
${\cal D}$ denotes the set of all non-increasing functions from $[0,1]$
to $\R$.

If $\sigma^2$ were known, then going from the white noise model (\ref
{eqmodel}) to the regression model (\ref{eqregmod}) would amount to
replace $F(t) = \int_0^t f(x) \,\d x $ and the rescaled Brownian motion
$n^{-1/2} W(t)$, $0 \le t \le1$, by the approximations $(1/n)\sum_{1
\le i \le j} f(i/n)$ and $(1/n) \sum_{1 \le i \le j}\epsilon_i$,
$1\le
j \le n$, respectively, so that a counterpart for $F_n$ in the
regression model is the continuous piecewise linear function
$F^{\mathrm{reg}}_n$ on $[0,1]$ interpolating between the points
\[
\biggl(\frac{j}{n}, \frac{1}{n} \sum_{0 \le i \le j}
Y_i \biggr),\qquad j=0,\ldots,n,
\]
where $Y_0=0$. Note that $F^{\mathrm{reg}}_n$ is nothing but the cumulative sum
diagram of the data $Y_i, 1\le i \le n$, with equal weights $w_i = 1/n$
(see also Barlow \textit{et al.} \cite{BBBB}).
Precisely, if $\sigma^2$ were known, then we would consider a finite
collection ${\cal C}_n$ of sub-intervals of $[0,1]$ and we would reject
$H_0$ if
\[
\max_{I\in{\cal C}_n} S^{\mathrm{reg},I}_n > r_{\alpha,n},
\]
where $r_{\alpha,n}$ is calibrated under the hypothesis that $f\equiv
0$ and where for all $I\in{\cal C}_n$,
%
\begin{equation}
\label{eqSnreg} S^{\mathrm{reg},I}_n= \sqrt{\frac{n}{\sigma^2|I|}}
\sup_{t\in I} \bigl(\widehat{F_n^{\mathrm{reg}}} ^I(t)-F^{\mathrm{reg}}_n(t)
\bigr).
\end{equation}
Since $\sigma^2$ is unknown, we need to estimate it. For ease of
notation, we assume that $n$ is even and we consider the estimator
\[
\hat\sigma^2=\frac{1}{n}\sum_{i=1}^{n/2}
(Y_{2i}-Y_{2i-1})^2.
\]
We introduce $\bar n=n/2$, $\sigma_0 =\sigma/\sqrt2$,
\[
\bar Y_i=(Y_{2i-1}+Y_{2i})/2 \quad\mbox{and}\quad \bar
\epsilon_i=(\epsilon_{2i-1}+\epsilon_{2i})/\sqrt2
\]
for all $i=1,\ldots,\bar n$, and we define on $[0,1]$ the function
\[
\bar f_n(t)=\bigl(f(t-1/n)+f(t)\bigr)/2,
\]
where $f(t)$ is defined in an arbitrary way for all $t\leq0$.
Thus, from the original model (\ref{eqregmod}), we deduce the
regression model
%
\begin{equation}
\label{eqregmod2} \bar Y_i = \bar f_n (i/\bar n ) +
\sigma_0\bar\epsilon_i,\qquad i=1,\ldots,\bar n,
\end{equation}
with the advantage that the observations $(\bar Y_i)_{1\leq i\leq\bar
n}$ are independent of $\hat\sigma^2$ (see the proof of Theorem \ref
{theolevelreg}). Note that $(\bar\epsilon_i)_{1\leq i\leq\bar n}$
are independent standard Gaussian variables, and a natural estimator
for $\sigma_0^2$ is $\hat\sigma^2_0=\hat\sigma^2/2$. Thus, the
regression model (\ref{eqregmod2}) has similar features as model
(\ref
{eqregmod}), so we proceed as described above, just replacing the
unknown variance by its estimator. Precisely, we choose some finite
collection ${\cal C}_{\bar n}$ of closed sub-intervals of $[0,1]$ with
endpoints on the grid $\{i/\bar n;i=0,\ldots,\bar n\}$ and for all
$I\in
{\cal C}_{\bar n}$, we define
\[
\hat S^{\mathrm{reg},I}_{\bar n} = \sqrt{\frac{\bar n}{\hat\sigma^2_0|I|}}
\sup_{t\in I} \bigl(\widehat{F_{\bar{n}}^{\mathrm{reg}}} ^I(t)-F^{\mathrm{reg}}_{\bar
n}(t)
\bigr),
\]
where $F^{\mathrm{reg}}_{\bar n}$ is the cumulative sum diagram of the data
$\bar Y_i, 1\le i \le\bar n$, with equal weights $w_i = 1/\bar n$. For
a given $\alpha\in(0,1)$, we reject $H_0\dvtx  f\in\cal D$ at level
$\alpha
$ when
%
\begin{equation}\label{eqregionreg}
\max_{I\in{\cal C}_{\bar n}} \hat S^{\mathrm{reg},I}_{\bar n} > r_{\alpha,n},
\end{equation}
where $r_{\alpha,n}$ is calibrated under the hypothesis that $f\equiv
0$. In order to describe more precisely $r_{\alpha,n}$, let us define
%
\begin{equation}
\label{Z} Z^{\mathrm{reg}}_n=\max_{I\in{\cal C}_{\bar n}} \sqrt{
\frac{\bar n}{|I|}} \sup_{t\in I} \bigl(\widehat{G}_{\bar
n}^I(t)-G_{\bar n}(t)
\bigr),
\end{equation}
where $G_{\bar n}$ is the cumulative sum diagram of the variables $\bar
\epsilon_i, 1\le i \le\bar n$, with equal weights $w_i = 1/\bar n$.
Although $G_{\bar n}$ is not observed, its distribution is entirely
known, and so is the distribution of $Z^{\mathrm{reg}}_n$.
We define $r_{\alpha,n}$ as the $(1-\alpha)$-quantile of
$Z^{\mathrm{reg}}_n/\sqrt{\chi^2(\bar n) /\bar n}$, where $\chi^2(\bar n)$ is a
random variable\vspace*{1pt} independent of $Z^{\mathrm{reg}}_n$ and having chi-square
distribution with $\bar n$ degrees of freedom. Approximated values for
the quantiles $r_{\alpha,n}$ can be obtained via Monte Carlo
simulations, and the test with critical region (\ref{eqregionreg}) is
of non-asymptotic level $\alpha$, as stated in the following theorem.
%
\begin{theo}\label{theolevelreg}
For every $\alpha\in(0,1)$,
\[
\sup_{f\in\cal D}\P_f \Bigl[\max_{I\in{\cal C}_{\bar n}} \hat
S_{\bar
n}^{\mathrm{reg},I}>r_{\alpha,n} \Bigr]=\alpha.
\]
Moreover, the above supremum is achieved at $f\equiv c$, for any fixed
$c\in\R$.
\end{theo}

As in Section \ref{secpower}, we study the performance of the test in
the case where
%
\begin{equation}
\label{eqCnreg} {\cal C}_{\bar n}= \biggl\{ \biggl[\frac{i}{\bar n},
\frac{j}{\bar
n} \biggr], i<j \mbox{ in }\{0,\ldots,\bar n\} \biggr\}.
\end{equation}
%
We obtain uniform separation rates that are comparable with the optimal
rates we have obtained in the white noise model, see Theorem \ref
{theoallpower} and Proposition \ref{proplowerbound}. For all $s\in
(1,2]$, $R>0$ and $L>0$,
denote
\[
{\cal F} (s,R,L) = {\cal F} (s,R) \cap\bigl\{f\dvtx [0,1] \rightarrow\R\mbox
{ s.t. }
\bigl\|f'\bigr\|_{\infty} \leq L\bigr\}.
\]

\begin{theo}\label{theopowerreg}
Let $\alpha,\beta$ in $(0,1)$, ${\cal C}_{\bar n}$ be the
collection (\ref{eqCnreg}) and $T^{\mathrm{reg}}_n$ be the test with critical
region~(\ref{eqregionreg}). Let $L>0$, $s\in(0,2]$ and $R>0$ and
assume $n\geq18\log(2/\alpha)$. In case $s\in(0,1]$, we assume that
$R/\sigma\leq n^s$ whereas in case $s\in(1,2]$, we assume that
$L/\sigma
\leq n$, ${\bar n}^s\sqrt{\log\bar n}\geq3^{s+1/2}R/\sigma$ and
$R/\sigma\geq2^{1+2s}\sqrt{(\log\bar n)/\bar n}$. Then, there
exists a
positive real $C(s,\alpha,\beta)$ only depending on $s$, $\alpha$ and
$\beta$ such that in case $s\in(0,1]$,
%
\begin{equation}
\label{eqpowerbhlreg} \rho\bigl(T^{\mathrm{reg}}_n,{\cal F}(s,R),
\Delta_1\bigr)\leq C(s,\alpha,\beta)R^{1/(1+2s)} \biggl(
\frac{\sigma^2\log n}{n} \biggr)^{s/(1+2s)}
\end{equation}
and, in case $s\in(1,2]$,
%
\begin{equation}
\label{eqpowerbhlreg2} \rho\bigl(T^{\mathrm{reg}}_n,{\cal
F}(s,R,L),\Delta_1\bigr)\leq C(s,\alpha,\beta)R^{1/(1+2s)}
\biggl(\frac{\sigma^2\log n}{n} \biggr)^{s/(1+2s)}
\end{equation}
and
%
\begin{equation}
\label{eqpowerdsreg} \rho\bigl(T^{\mathrm{reg}}_n,{\cal F}(s,R,L),
\Delta_2\bigr)\leq C(s,\alpha,\beta)R^{3/(1+2s)} \biggl(
\frac{\sigma^2\log n}{n} \biggr)^{(s-1)/(1+2s)}.
\end{equation}
\end{theo}

\section{Simulation study and power comparison}\label{secsimulation}

In this section, we are first concerned with some algorithmic aspects
as well as with power comparisons. First, we explain how to compute
approximate values of our test statistics
\[
S_n:=\max_{I\in\mathcal C_n} S_n^I
\quad\mbox{or}\quad S^{\mathrm{reg}}_{n}:=\max_{I\in\mathcal C_{\bar n}} \hat
S^{\mathrm{reg},I}_{\bar n}.
\]
We give also a brief description of how the critical thresholds and the
power are calculated. For $n=100$ we study the power of our test under
various alternatives in the white noise and Gaussian regression models
so as to provide a comparison with other monotonicity tests such as
those proposed by Gijbels \textit{et al.} \cite{gijbelskoch00} and Baraud
\textit{et al.}
\cite{bhl05}. Although the tests are expected to behave better for
large sample sizes, we give below power results in the regression model
for $n=50$ and some of the examples considered in the aforementioned
papers to give an idea of the performance of the corresponding test for
moderate sample sizes.

\subsection{Implementing the test}

\subsubsection*{Computing $S_n$} We describe here the numerical procedures
used to compute our statistic for the white noise model. Some of these
numerical procedures take a simpler form for the Gaussian regression
model, and hence they will be only briefly described below. In the
white noise model, the observed process $F_n$ can be only computed on a
discrete grid of $[0,1]$. For each subinterval $I$ of $[0,1]$, let
$\tilde
{F}^I_n$ be the approximation of the restriction of $F_n$ to $I$
obtained via a linear interpolation between the values of $F_n$ at the
points of the chosen grid. Also, let $\widehat{\tilde{F}}{}^I_n$ be its
least concave majorant. For a given pair of integers $(i, j)$ such that
$0 \le i < j
\le n$, we use the Pool Adjacent Violators Algorithm (see, e.g., \cite
{BBBB}) to compute the slopes and vertices of $\widehat{\tilde
{F}}{}^{I_{ij}}_n$, where $I_{ij} = [i/n, j/n]$. In order to gain in
numerical efficiency, we compute the least concave majorants
progressively, taking
advantage of previous computations\vspace*{1pt} on smaller intervals using a
concatenation procedure. More precisely, we first compute $\widehat
{\tilde{F}}{}^{I_{j(j+1)}}_n$, for all $j \in\{0,\ldots,n-1\}$, and store
the corresponding sets of vertices. Then, for $l \in\{j+2,\ldots, n \}$,
the least concave majorant $\widehat{\tilde{F}}{}^{I_{jl}}_n$ is also
equal to the
least concave majorant
of the function resulting from the concatenation of $\widehat{\tilde
{F}}{}^{I_{j(l-1)}}_n$ and $\widehat{\tilde{F}}{}^{I_{(l-1)l}}_n$, whose
sets of vertices were previously stored. This merging step reduces the
computation time substantially, because the number of vertices of a
least concave majorant on $I_{ij} $ is often much smaller than the
number of grid
points in $I_{ij}$. At each step of the algorithm, the maximum of
$\widehat{\tilde{F}}{}^{I_{jl}}_n - \tilde{F}^{I_{jl}}_n$ on $I_{jl}$ is
multiplied by $n/\sqrt{(l-j)}$ and stored. Last, we obtain the
(approximated) value of the test statistic $S_n$ by taking the largest
value of those rescaled maxima and then dividing by $\sigma$.

\subsubsection*{Computing the statistic $S^{\mathrm{reg}}_{n}$} The discretized
nature of this setting makes the computations faster than in the white
noise model. The least concave majorants based on the independent data
$\bar{Y}_i, i=1,\ldots, \bar n$ are progressively computed on $I_{ij} = \{i/\bar n,\ldots, j/\bar n\}, 0\le i < j \le\bar n $, using the concatenation
technique as described above. Note that each data point $\bar{Y}_i$ is
assigned to the design point $x_i = i/\bar n, i=1,\ldots, \bar n$, so
that only half of the original grid is exploited. This is the price to
be paid for not knowing the variance of the noise. The maximum
deviation between the cumulative sum diagrams and their corresponding
least concave majorants yields the value of
the test statistic $S^{\mathrm{reg}}_n$ after division by the estimate $\widehat
\sigma_0 = (\bar n^{-1} \sum_{j=1}^{\bar n} (Y_{2i -1} - Y_{2i})
)^{1/2}$.

\begin{figure}

\includegraphics{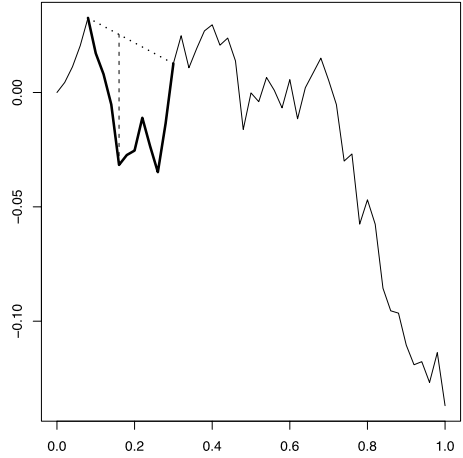}

\caption{
The plot of the cumulative sum diagram of $\bar Y_i, i=1,\ldots, \bar
n=50$ based on 100 independent realizations of standard Gaussians $Y_1,\ldots, Y_{100}$, and the least concave majorant on $I= \{4/50,\ldots,
15/50 \}$
yielding the maximal value of deviation $Z^{\mathrm{reg}}_n$.}
\label{Zreg}
\end{figure}

In Figure \ref{Zreg}, we illustrate the computation of $Z^{\mathrm{reg}}_n$.
Independent replications of the above calculations under the hypothesis
$f\equiv0$ enable us to compute the empirical quantiles of $S^{\mathrm{reg}}_n$.
For a given $\alpha\in(0,1)$, the empirical quantile of order
$1-\alpha$ will be taken as an approximation for the critical threshold
for the asymptotic level $\alpha$, which will be denoted by $r_{\alpha,n}$.

\subsubsection*{Computing the critical thresholds and the power} 
\begin{figure}

\includegraphics{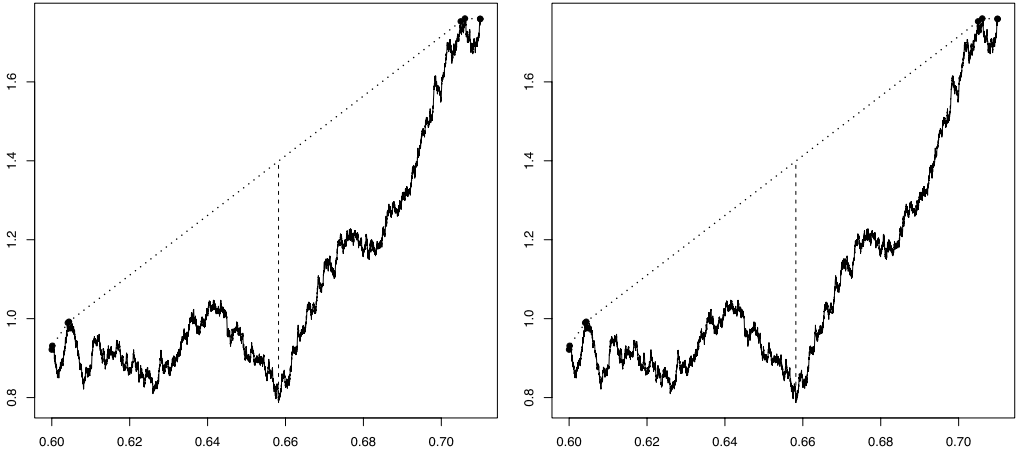}

\caption{Left plot: Brownian approximation on the unit interval and the
least concave majorant on the subinterval $I = [61/100, 72/100]$
yielding the
approximate value of $S_n$ for $f\equiv0$. Right plot: magnified plot
of the least concave majorant on $I$. The vertices are shown in
bullets.}
\label{BrownMotion}
\end{figure}

We now
describe how we determine the critical region of our tests for a given
level $\alpha\in(0,1)$. The calculation of the power under the
alternative hypothesis is performed along the same lines, hence its
details are skipped. Determining the critical region relies on
computing Monte Carlo estimates of $s_{\alpha, n}$ and $r_{\alpha,n}$,
the $(1-\alpha)$-quantiles of the statistic $S_n$ and $S^{\mathrm{reg}}_n$ under
the least favorable hypothesis $f\equiv0$.
The calculations in the regression setting have been described above
and are much simpler than for the white noise model. Thus, we only
provide some details of how the approximation of the critical threshold
$s_{\alpha,n}$ is performed. Approximation of $s_{\alpha,n}$ requires
simulation of $C$ independent copies of Brownian motion on $[0,1]$. For
a chosen $r\in\mathbb N^\star$, we simulate $m=n\times r$ independent
standard Gaussians $Y_1,\ldots,Y_m$. The rescaled partial sums
$m^{-1/2}\sum_{i=1}^k Y_i$, for $k=0,\ldots,m$, provide approximate
values for Brownian motion at the points of the regular grid $\{k/m;
k=0,\ldots,m\}$ of $[0,1]$. We then proceed as explained in the
previous paragraph to obtain the approximations $\tilde W^{I_{ij}}$,
$0\leq i <j\leq n$, of the restrictions of Brownian motion to the
intervals $I_{ij}$ and their respective least concave majorants
$\widehat{\tilde{W}}{}^{I_{ij}}$.
In the sequel, we fix $n=100$. Figure \ref{BrownMotion} shows an
example where approximated value of $S_n$ for $f\equiv0$ is found to
be equal to $ 0.611 \times\sqrt{100/11} \approx1.842$, where $0.611$
is the length of the vertical dashed line representing the maximal
difference between $\tilde{W}^I$ and its least concave majorant. For
the approximation
of Brownian motion, we have taken $r=1000$.

\begin{table}[b]
\tablewidth=265pt
\caption{Monte Carlo estimates of $s_{\alpha, n}$, the $(1-\alpha
)$-quantiles of $S_n$ for $f\equiv0$. The estimation is based on
$C=5000$ runs, $n=100$ and $r=1000$}
\label{flottantsQuantilestable}
\begin{tabular*}{\tablewidth}{@{\extracolsep{\fill}}llllll@{}}
\hline
$\alpha$ & 0.01 & 0.02 & 0.03 & 0.04 & 0.05 \\
$s_{\alpha,n}$ & 2.451860 & 2.384279 & 2.343395 & 2.308095 &
\textbf{2.278482} \\
[6pt]
$\alpha$ & 0.06 & 0.07 & 0.08 & 0.09 & 0.10 \\
$s_{\alpha,n}$ & 2.254443 & 2.235377 & 2.217335 & 2.205137 & 2.191502
\\
\hline
\end{tabular*}
\end{table}

Based on $C=5000$ runs, we found that $r=1000$ and $10\,000$ yield
close results for the distribution of the approximated value of the
test statistic $S_n$. On the other hand, larger values of $r$ make the
computations prohibitively slow. Thus, we chose $r=1000$ as a good
compromise for Monte Carlo estimation of the quantiles as well as for
power calculations. Based on $C=5000$ runs, Monte Carlo estimates of
$s_{\alpha, n}$, the $(1-\alpha)$-quantiles of $S_n$ for $f\equiv0$,
were computed for $\alpha\in\{0.01, 0.02, \ldots, 0.1\}$ and are
gathered in Table \ref{flottantsQuantilestable}. In particular, the
Monte Carlo estimate of the quantile of order $1-\alpha= 0.95$ is
found to be equal to 2.278482. Note that the true quantile should be
comprised between $q(\alpha)$ and $q (2\alpha[n(n+1)]^{-1}
)$, the $1-\alpha$ and $1 - 2\alpha[n(n+1)]^{-1}$ quantiles of $Z$
(cf. Section~\ref{seclevel}). A numerical method for finding very
precise approximations of upper quantiles of $Z$ was developed by
Balabdaoui and Filali \cite{balabdfilali10} using a Gaver--Stehfest
algorithm. For $\alpha= 0.05$, the Gaver--Stehfest approximation of
$q(\alpha)$ and
$q (2\alpha[n(n+1)]^{-1} )$ computed yields the values
1.46279052 and 2.60451660, respectively. Hence, the obtained Monte
Carlo estimate seems to be consistent with the theory.

\begin{table}
\tablewidth=265pt
\caption{Monte Carlo estimates of $r_{\alpha, n}$, the $(1-\alpha
)$-quantiles of $S^{\mathrm{reg}}_n$ for $f\equiv0$. The estimation is based on
$C=5000$ runs and $n=100$}
\label{flottantsQuantilesregtable2}
\begin{tabular*}{\tablewidth}{@{\extracolsep{\fill}}llllll@{}}
\hline
$\alpha$ & 0.01 & 0.02 & 0.03 & 0.04 & 0.05 \\
$r_{\alpha,n}$ & 2.150903 & 2.090423 & 2.013185 & 1.998276 &
\textbf{1.970304} \\
[6pt]
$\alpha$ & 0.06 & 0.07 & 0.08 & 0.09 & 0.10 \\
$r_{\alpha,n}$ & 1.950475 &1.938510 & 1.906807 & 1.892049 & 1.870080 \\
\hline
\end{tabular*}
\end{table}

For testing monotonicity in the regression model in (\ref{eqregmod2}),
Table \ref{flottantsQuantilesregtable2} give values of $r_{\alpha, n}$
for $n=100$ and $\alpha\in\{0.01, 0.02,\ldots, 0.1 \}$. The
approximated quantiles are obtained with $C=5000$ independent draws
from a standard Gaussian. The Monte Carlo estimate of the quantile of
order $1-\alpha= 0.95$ is found to be equal to $r_{0.05,100}
=1.970304$, and hence smaller than its counterpart in the white noise
model. Recall that the statistic $S^{\mathrm{reg}}_n$ has the same
distribution as $Z^{\mathrm{reg}}_n/\sqrt{\chi^2(\bar n)/\bar n}$.
Thus, as $n$ gets larger, we expect the distribution of
$S^{\mathrm{reg}}_n$ to be closer to that of~$S_n$. We have added Table
\ref{flottantsQuantilesregtable} where we give the obtained values of
$r_{\alpha, n}$ for $n=1000$, and which are clearly close to the
approximated quantiles obtained in Table \ref {flottantsQuantilestable}
in the white noise model.


\subsection{Power study}

In this subsection, we shall determine the power of our tests when the
true signal deviates either globally or locally from monotonicity. The
functions we consider here have already been used by Gijbels, Hall,
Jones and Koch \cite{gijbelskoch00} and Baraud, Huet and Laurent \cite
{bhl05}. Our goal is then two-fold: compute the power of our test, and
compare its performance to that of the tests considered in these papers.

\begin{table}
\tablewidth=265pt
\caption{Monte Carlo estimates of $r_{\alpha, n}$, the $(1-\alpha
)$-quantiles of $S^{\mathrm{reg}}_n$ for $f\equiv0$. The estimation is based on
$C=5000$ runs and $n=1000$}
\label{flottantsQuantilesregtable}
\begin{tabular*}{\tablewidth}{@{\extracolsep{\fill}}llllll@{}}
\hline
$\alpha$ & 0.01 & 0.02 & 0.03 & 0.04 & 0.05 \\
$r_{\alpha,n}$ & 2.475129 & 2.405715 & 2.351605 & 2.320332 & 2.295896
\\[6pt]
$\alpha$ & 0.06 & 0.07 & 0.08 & 0.09 & 0.10 \\
$r_{\alpha,n}$ & 2.277919 &2.249911& 2.232569 & 2.214476 & 2.202449\\
\hline
\end{tabular*}
\end{table}

\begin{table}
\tablewidth=\textwidth
\caption{Power of the tests based on $S_n$, $S^{\mathrm{reg}}_n$ and $T_B$ for
$n=100$ (see text for details)}
\label{flottantsPowBaraudtable}
\begin{tabular*}{\tablewidth}{@{\extracolsep{\fill}}lllll@{}}
\hline
Function & $\sigma^2$ & $S_n$ & $S^{\mathrm{reg}}_n$ & $T_B$ \\
\hline
$f_1$ & 0.01& 1.00 & 0.99 & 0.99\\
$f_2$ & 0.01 & 0.99 & 1.00 & 0.99\\
$f_3$ & 0.01& 1.00 & 0.98 & 1.00 \\
$f_4$ & 0.01 & 0.99 & 0.99 & 0.99\\
$f_5$ & 0.004 & 1.00 & 0.99 & 0.99\\
$f_6$ & 0.006 & 1.00 & 0.99 & 0.98 \\
$f_7 $ & 0.01 & 0.79 & 0.68 & 0.76 \\
\hline
\end{tabular*}
\end{table}

\begin{table}[b]
\caption{Power of the tests based on $S_n$, $S^{\mathrm{reg}}_n$ and
$T_{\mathrm{run}}$ for $n=100$ (see text for details)}
\label{flottantsPowGijbelstable}
\begin{tabular*}{\tablewidth}{@{\extracolsep{\fill}}lllllllllll@{}}
\hline
& &\multicolumn{3}{l}{$a=0$} & \multicolumn{3}{l}{$a=0.25$} &
\multicolumn{3}{l@{}}{$a=0.45$} \\[-4pt]
& &\multicolumn{3}{l}{\hrulefill} & \multicolumn{3}{l}{\hrulefill} &
\multicolumn{3}{l@{}}{\hrulefill} \\
& $\sigma$ & 0.025 & 0.05 & 0.1 & 0.025 & 0.05 & 0.1 & 0.025 & 0.05 &
0.1 \\
\hline
$S_n$ & & 0.010 & 0.018 & 0.014 & 0.246 & 0.043 & 0.031 & 1.000 & 1.000
& 0.796 \\
$S^{\mathrm{reg}}_n$ & & 0.000 & 0.002& 0.013 & 0.404 & 0.053 & 0.007& 1.000
&1.000 & 0.683 \\
$T_{\mathrm{run}}$ & & 0.000 & 0.000 & 0.000 & 0.106 & 0.037 & 0.014 &
1.000 & 1.000 & 0.805 \\
\hline
\end{tabular*}
\end{table}

In Tables \ref{flottantsPowBaraudtable} and \ref
{flottantsPowGijbelstable} below, the columns labeled $S_n$ and
$S^{\mathrm{reg}}_n$ give the power of the tests of level $\alpha= 0.05$ based
on the statistics $S_n$ and $S^{\mathrm{reg}}_n$ in the white noise and
regression model (\ref{eqregmod}), respectively, that is, the proportion
of $S_n > s_{\alpha, n}$ and $S^{\mathrm{reg}}_n > r_{\alpha, n}$ among the
total number of runs.

\subsubsection*{\texorpdfstring{Functions considered by Baraud
\textit{et al.} \cite{bhl05}}{Functions considered by Baraud \textit{et al.} [5]}}

Baraud, Huet and Laurent \cite{bhl05} consider a regression model with
deterministic design points $0 \le x_1 \le\cdots\le x_n \le1$ and
Gaussian noise whose variance is finite and equal to $\sigma^2$. Their
monotonicity test of the true regression function is based on
partitioning the set $\{1,\ldots, n\}$ into $l \in\{2,\ldots, l_n \}$
subintervals for a given integer $2 \le l_n \le n$. Below, we use $T_B$
as a shorthand notation for their local mean test where the maximal
number of subsets in the partition is $l_n = 25$ and $x_i = i/n$,
$i=1,\ldots, n$. The basis of power comparison consists of the
following functions
\begin{eqnarray*}
f_1(x) &= & -15 (x-0.5)^3 1_{x \le0.5} - 0.3
(x-0.5) + \exp\bigl(-250(x-0.25)^2\bigr),
\\
f_2(x) & = & 1.5 \sigma x,
\\
f_3(x) & = & 0.2 \exp\bigl(-50(x-0.5)^2\bigr),
\\
f_4(x) & = & -0.1 \cos(6 \uppi  x),
\\
f_5(x) &= & -0.2 x + f_3(x),
\\
f_6(x) & = & -0.2 x + f_4(x),
\\
f_7(x) & = & -(1+x) + 0.25 \exp\bigl(- 50 (x-0.5)^2
\bigr).
\end{eqnarray*}
Note that $f_7$ is a special case of Model III of Gijbels \textit{et al.}
\cite
{gijbelskoch00}, which we consider below.


Table \ref{flottantsPowBaraudtable} gathers power results for the
functions $f_1/\sigma_1,\ldots, f_7/\sigma_7$ where $\sigma^2_i$ is the
variance of the noise considered by the authors when the true function
is $f_i$ (see the second column in Table \ref
{flottantsPowBaraudtable}). Calculation of the power of our local
least concave majorant test was based on 1000 independent runs in both
white noise and
regression models. We see that in both models our tests perform as well
as the local mean test of Baraud \textit{et al.} \cite{bhl05}, except
for the
function $f_7$, which is a special case of Model III considered by
Gijbels \textit{et al.} \cite{gijbelskoch00} ($a = 0.45$ and $\sigma=0.1$),
and where our test in the regression model seems to be doing a bit
worse than local mean test. We would like to note that in the white
noise model the power obtained for the functions $f_1,\ldots, f_6$ is not
much altered when we replace the quantile $s_{\alpha, n}$ by the
upper-bound $q (2\alpha[n(n+1)]^{-1} ) \approx2.60451660$
(see Lemma \ref{lemmajq}): in this case, the power is slightly smaller
and we find a minimum difference of order $-$0.070.

\subsubsection*{\texorpdfstring{Functions considered by Gijbels \textit{et al.}
\cite{gijbelskoch00}}{Functions considered by Gijbels \textit{et al.}
[16]}}

Gijbels \textit{et al.} \cite{gijbelskoch00} consider a regression
model where
the deviation from monotonicity of the true regression function $f$
defined on $[0,1]$ depends on a parameter $a > 0$ in the following manner
\[
f_a(x) = - (1 + x) + a \exp\bigl(- 50 (x-0.5)^2\bigr),\qquad
x \in[0,1].
\]
(Note that we have multiplied their function by $(-1)$ to obtain a
perturbation of a decreasing function, here $-1 - x$.) We refer to
\cite
{gijbelskoch00} for a description of their tests.
Here, we compare the performance of our test to their test based on the
statistic $T_{\mathrm{run}}$ (runs of equal signs).

For $a=0$, where the true function $-1-x$ satisfies $H_0$, all three
tests have a rejection probability which is much smaller than $\alpha
=0.05$. Our test in the white noise model seems to be however less
conservative than our test in the regression model and the test of
Gijbels \textit{et al.} \cite{gijbelskoch00}. Our power values in both
models (\ref{eqmodel}) and (\ref{eqregmod}) suggest that our tests
are exhibiting comparable performance, except for the configuration
$a=0.45$ and $\sigma=0.1$ where our test based on $S_n^{\mathrm{reg}}$ seems to
be doing a bit worse. This fact was also noted above in our comparison
with the local mean test of Baraud \textit{et al.} \cite{bhl05}. In the white
noise model, we would like to note that replacing the quantile $s_{\alpha,
n}$ by $q (2\alpha[n(n+1)]^{-1} )$ implies now a strong
decrease in the power. We conclude that the upper-bound of Lemma \ref
{lemmajq}, although interesting in its own right, should be used with caution.

\begin{table}
\tablewidth=165pt
\caption{Power of the test based on $S^{\mathrm{reg}}_n$ for $n=50$ (see text
for details)}\label{PowBaraudn50}
\begin{tabular*}{\tablewidth}{@{\extracolsep{\fill}}lll@{}}
\hline
Function & $\sigma^2$ & Power \\
\hline
$f_1$ & 0.01& 0.44 \\
$f_2$ & 0.01 & 0.68 \\
$f_3$ & 0.01& 0.84 \\
$f_4$ & 0.01 & 0.36 \\
$f_5$ & 0.004 & 0.90 \\
$f_6$ & 0.006 & 0.44 \\
$f_7 $ & 0.01 & 0.34 \\
\hline
\end{tabular*}
\end{table}

Finally, we have carried out further investigation of the performance
of our test in the regression model for the moderate sample size
$n=50$. The Monte Carlo estimate of the 95\% quantile of $S^{\mathrm{reg}}_n$ is
found to be approximately equal to $1.837931$. Table \ref{PowBaraudn50}
gives the obtained values of the power for the alternatives considered
by Baraud \textit{et al.} \cite{bhl05}. As expected, the test loses
from its
performance for this smaller sample size but is still able to detect
deviation from monotonicity for some of those alternatives with large
power (functions $f_3$ and $f_5$). The alternatives considered by
Gijbels \textit{et al.}  \cite{gijbelskoch00} present a much bigger challenge and
the power of our test is found to be small with values between 0.003
and 0.01 for $a=0.25$. For $a=0.45$, our test is found to be powerful
with power values equal to 1 and 0.954 for $\sigma= 0.025$ and $0.05$,
respectively.


\section{Discussion}\label{secdiscussion}
In this section, we compare our test method with some of its (existing)
competitors and also discuss possible generalizations to other settings.

\subsection{Comparison with competing tests}

In this article, we have proposed a new procedure for testing
monotonicity which is able to detect localized departures from
monotonicity, with exact prescribed level, either in the white noise
model or in the Gaussian regression framework with unknown variance.
Firstly, as explained in Section \ref{secsimulation}, our test
statistic, which is a maximum over a finite number of intervals --
increasing with $n$ -- of local least concave majorants, can be
computed exactly and efficiently. This is a big advantage when our
method is compared with the two test statistics proposed in \cite
{ds01}, as they rely on a family of kernel estimators indexed by a
continuous set of bandwidths which must be in practice discretized.
Secondly, the distribution of our test statistic under the least
favorable hypothesis is known and its quantiles can be evaluated via
Monte Carlo simulations. We insist on the fact that, as opposed
to \cite
{bhl05}, only one quantile has to be evaluated for a given level, and
that no bootstrap is required as opposed to \cite{hallheckman00}.
Moreover, our test statistic does not rely on any smoothness
assumption, because no smoothing parameter is involved in the
construction of the test as opposed to \cite{ghosal00}.
In terms of power, an interesting property of our procedure is its
adaptivity. Indeed, our procedure attains the same rates of separation
as \cite{ghosal00} (without having to play with some smoothness
parameter), as well as the optimal rates obtained in the four
procedures in \cite{bhl05,ds01}. Lastly, the detailed comparative study
above shows that the power values we attain are generally similar to
those obtained by \cite{bhl05} and better than the ones obtained by
\cite{gijbelskoch00}. Note also that in practice our procedure reaches
at most the prescribed level, which may not be the case for \cite
{hallheckman00}, even for Gaussian errors. We conclude that our
approach seems to enjoys all the qualities of the existing methods.

\subsection{About the choice of ${\cal C}_{n}$}

As explained at the beginning of Section \ref{secpower}, our choice of
${\cal C}_{n}$ is motivated by our wish to detect as many alternatives
as possible. However, if one knows in advance properties of the
subinterval where $f$ is likely to violate the non-increasing
assumption, then one can incorporate this knowledge to the choice of
${\cal C}_{n}$, that is, one can choose a reduced collection of
subintervals in such a way that for the largest interval $I\subset
[0,1]$ where $f$ is likely to violate the non-increasing assumption,
there is an interval close to $I$ in the chosen collection ${\cal
C}_{n}$. By ``reduced collection'', we mean a collection that is
included in the one defined in (\ref{eqCn}), and by ``an interval
close to $I$'', we mean, for instance, an interval whose intersection
with $I$ has a length of the same order of magnitude (up to a
multiplicative constant that does not depend on $n$) as the length of
$I$, and where the increment of $f$ has the same order as on $I$. It
can be seen from our proofs that for an arbitrary choice of ${\cal
C}_{n}$, we obtain that there exists $C(\alpha,\beta)>0$ only depending
on $\alpha$ and $\beta$ such that (\ref{eqpower}) holds provided there
exists $[x,y]\in{\cal C}_{n}$ such that
%
\begin{equation}
\label{eqpowerslopebis}
\sup_{t\in[x,y]}\frac{t-x}{\sqrt{y-x}}(\bar
f_{xy}-\bar f_{xt})\geq C(\alpha,\beta)\sqrt
\frac{{\sigma^2\log}{|\cal C}_{n}|}{n}.
\end{equation}
Note that in the case where ${\cal C}_{n}$ is given by (\ref{eqCn}),
$|{\cal C}_{n}|=n(n+1)/2$ so that $\log|{\cal C}_{n}|$ is of the order
$\log n$, hence the $\log n$ term in Theorem \ref{theopowerslope}. In
view of condition (\ref{eqpowerslopebis}), good power properties are
obtained if one chooses ${\cal C}_{n}$ in such a way that $|{\cal
C}_{n}|$ is not too large, but ${\cal C}_{n}$ contains an interval
$[x,y]$ close to the largest interval where $f$ is likely to violate
the non-increasing assumption: ${\cal C}_{n}$ must have good
approximation properties while having a moderate cardinality.
For instance, to test that $f$ is non-increasing on $[0,1]$ against the
alternative that $f$ is U-shaped on $[0,1]$, one can consider a
collection of intervals of the form $[x,1]$. However, in such a case,
only alternatives that are increasing on the right boundary of $[0,1]$
could be detected. More generally, considering a reduced collection
${\cal C}_{n}$ may cause a loss of adaptivity so we do not pursue the
study of our test in the case where it is defined with a reduced
collection ${\cal C}_{n}$.

\subsection{Possible extensions to more general models}
Recall that in the case where we observe $Y_{1},\ldots,Y_{n}$ according
to (\ref{eqregmod}) with a known $\sigma>0$ and i.i.d. standard
Gaussian $\epsilon_{i}$'s, we reject $H_{0}\dvtx  f\in{\cal D}$ against the
alternative $H_1\dvtx  f\notin{\cal D}$ if
%
\begin{equation}
\label{eqrejet}
\max_{I\in{\cal C}_{n}}S_{n}^{\mathrm{reg},I}>r_{\alpha,n},
\end{equation}
where $S_{n}^{\mathrm{reg},I}$ is given by (\ref{eqSnreg}) and $r_{\alpha,n}$
is the $(1-\alpha)$-quantile of
$\max_{I\in{\cal C}_{n}}S_{n}^{\mathrm{reg},I}$ under the hypothesis that
$f\equiv0$.

The hypothesis that the $\epsilon_{i}$'s are standard Gaussian can
easily be relaxed provided that the common distribution remains known.
Indeed, provided both $\sigma$ and the common distribution of the
$\epsilon_{i}$'s are known, the law of $\max_{I\in{\cal
C}_{n}}S_{n}^{\mathrm{reg},I}$ is entirely known (at least theoretically) under
the least favorable hypothesis that $f\equiv0$, so we can obtain
approximate value of the $(1-\alpha)$-quantile $r_{\alpha,n}$ of this
law via Monte Carlo simulations. The critical region (\ref{eqrejet})
then defines a test with level~$\alpha$. If there exists $C>0$ and
$C'>0$ such that for all $x>0$ and $I=[i/n,j/n]\in{\cal C}_{n}$,
\[
\P\Biggl(\max_{i<k\leq j}\sum_{l=i+1}^k
\delta_{i}>x \sqrt{j-i} \Biggr)\leq C\exp\bigl(-C'x^2
\bigr),
\]
where $\delta_i$ denotes either $\epsilon_{i}$ or $-\epsilon_{i}$, then
this test has similar power properties as in the Gaussian case (the
rates in Theorem \ref{theopowerreg} still hold, with possibly
different constants).

In the case where the distribution of the $\epsilon_{i}$'s is known but
$\sigma$ is unknown, it is tempting to replace $\sigma$ with an
estimator $\hat\sigma_{n}$ in the definition of the test statistic, and
to consider the critical region
$T_{n}/\hat\sigma_{n}>\hat r_{\alpha,n}$
where
\[
T_n=\sigma\max_{I\in{\cal C}_{n}}S_{n}^{\mathrm{reg},I}
\]
is observable and $\hat r_{\alpha,n}$ is the $(1-\alpha)$-quantile of
$T_{n}/\hat\sigma_{n}$
under the hypothesis that $f\equiv0$. However, the distribution of
$T_{n}/\hat\sigma_{n}$ is not known in the general case even if
$f\equiv
0$. In the particular case where the $\epsilon_{i}$'s are standard
Gaussian, we suggest in Section \ref{secregression} to rather consider
\[
\frac{\sigma_{0}}{\hat\sigma_{0}}\max_{I\in{\cal C}_{\bar
n}}S_{\bar n}^{\mathrm{reg},I}
\]
as a test statistic, where $\hat\sigma_{0}$ and $S_{\bar n}^{\mathrm{reg},I}$
are defined in such a way that the distribution of $\hat\sigma
_{0}/\sigma_{0}$ is known if $f\equiv0$, and $\hat\sigma_{0}$ is
independent of $\max_{I\in{\cal C}_{\bar n}}S_{\bar n}^{\mathrm{reg},I}$. This
way, the distribution of the test statistic is known under the
hypothesis that $f\equiv0$ and we are able to calibrate the test.
Except in the Gaussian case, the distribution of the test statistic is
not known even under the hypothesis that $f\equiv0$. In such
situations, it is tempting to argue asymptotically, as $n\to\infty$.
This requires computation of the limit distribution of the test
statistic $T_{n}/\hat\sigma_{n}$ under the least favorable hypothesis.
More precisely, suppose there exist sequences $a_n$ and $b_{n}$ such
that if $f\equiv0$,
\[
a_{n}(T_{n}/\sigma-b_{n})
\]
converges in distribution to $T$ as $n\to\infty$, where $T$ is a random
variable with a continuous distribution function. Suppose, moreover,
that $\hat\sigma_n$ converges in probability to $\sigma$ as $n\to
\infty
$, and either $a_{n}b_{n}=\RMO(1)$ or $\hat\sigma_{n}=\sigma
+\RMo_P(1/(a_{n}b_{n}))$. Denoting by $s_\alpha$ the $(1-\alpha)$-quantile
of $T$, the test defined by the critical region
%
\begin{equation}
\label{eqrejetextension}
T_{n}/\hat\sigma_{n}>b_{n}+a_n^{-1}s_{\alpha}
\end{equation}
has asymptotic level $\alpha$. This means that one can extend our
method (even in the case where the distribution of the $\epsilon_{i}$'s
is unknown) by considering a critical region of the form~(\ref
{eqrejetextension}).

Other extensions of the method are conceivable to test, for instance,
that a density function or a failure rate is non-increasing on a given
interval, or that a regression function is non-increasing in a
heteroscedastic model. In such cases, the exact distribution of the
test statistic is not known even under the least favorable hypothesis
since it depends on unknown parameters. Similar to the regression case
with non-Gaussian errors above, it is then tempting to replace the
unknown parameters with consistent estimators and argue asymptotically.
Such asymptotic arguments, with the computation of the limit
distribution of the test statistic, is beyond the scope of the present
paper and is left for future research.

\section{Proofs}\label{secallproofs}

Without loss of generality (see \cite{CFNsupp}), we assume for
simplicity that the noise level is $\sigma=1$.

\subsection{\texorpdfstring{Proofs for Section \protect\ref{seclevel}}
{Proofs for Section 2}}\label
{secproofslevel}

\begin{pf*}{Proof of Lemma \ref{lemleastfav}}
Define
$B(t)= (W(t(b-a)+a)-W(a) )/\sqrt{|I|}$ on $[0,1]$.
From Lem\-ma~2.1 in \cite{dt03}, it follows that for all $t\in I$,
\[
\hat W^I(t)-W(t) =\sqrt{|I|} \biggl[\hat B \biggl(\frac{t-a}{b-a}
\biggr)-B \biggl(\frac
{t-a}{b-a} \biggr) \biggr].
\]
If $f\equiv c$ on $I$ then $F$ is linear on $I$, so Lemma 2.1 in \cite
{dt03} shows that
$\hat F_n^I=F+\hat W^I/\sqrt n$ and
\[
S_n^I= \frac{1}{\sqrt{|I|}}
\sup_{t\in I} \bigl(\hat W^I(t)-W(t) \bigr)=\sup_{t\in[0,1]}
\bigl(\hat B(t)-B(t) \bigr).
\]
But $B$ is distributed like $W$, so $S_n^I$ is distributed like $Z$.
\end{pf*}
\begin{pf*}{Proof of Theorem \ref{theoIlevel}}
For every $f\in{\cal D}^I$, $F$ is concave on $I$, so the process
$F+\hat W^I/\sqrt n$ is concave and above $F_n$ on $I$. Thus, it is
also above $\hat F_n^I$, and hence
\[
S_n^I \leq\sqrt\frac{n}{|I|}\sup_{t\in I}
\biggl(F(t)+\frac{1}{\sqrt
n}\hat W^I(t)-F_n(t) \biggr)
\leq\frac{1}{\sqrt{|I|}}\sup_{t\in I} \bigl(\hat W^I(t)-W(t)
\bigr).
\]
Since $Z$ is continuously distributed (see \cite{durot03}, Lemma 1),
Lemma \ref{lemleastfav} yields
\[
\sup_{f\in{\cal D}^I} \P_f \bigl[S_n^I>q(
\alpha) \bigr]\leq\P\bigl[Z>q(\alpha) \bigr]=\alpha.
\]
Now, suppose $f\equiv c$ over $I$ for a fixed $c\in\R$. Then Lemma
\ref
{lemleastfav} shows that
\[
\P_f \bigl[S_n^I>q(\alpha) \bigr]=\P
\bigl[Z>q(\alpha) \bigr]=\alpha.
\]
\upqed
\end{pf*}
%
%
\begin{pf*}{Proof of Theorem \ref{theolevel}}
Suppose $f\in\cal D$. We deduce as in the proof of Theorem \ref
{theoIlevel} that
\[
\max_{I\in{\cal C}_n}S_n^I\leq\max_{I\in{\cal C}_n}\sqrt{
\frac
{1}{|I|}}\sup_{t\in I}\bigl(\hat W^I(t)-W(t)\bigr)
\]
with equality when $f\equiv c$ for some fixed $c\in\R$, hence Theorem
\ref{theolevel}.
\end{pf*}
\begin{pf*}{Proof of Lemma \ref{lemmajq}}
It follows from Lemma \ref{lemleastfav} that for every $x\geq0$,
\[
\P\biggl[\max_{I\in{\cal C}_n}\sqrt{\frac{1}{|I|}}\sup_{t\in
I}\bigl(
\hat W^I(t)-W(t)\bigr)>x \biggr]\leq|{\cal C}_n|\P[Z>x
].
\]
Thus, with $x=q ({\alpha}/|{\cal C}_n| )$, we obtain the first
inequality.
Now, from the definition (\ref{eqZ}) of $Z$, we have
\[
Z\leq\sup_{t\in[0,1]} W(t)+\sup_{t\in[0,1]}\bigl(-W(t)\bigr).
\]
Since $W$ has the same distribution as $-W$ and its supremum over
$[0,1]$ satisfies an exponential inequality (see, e.g., \cite{ry91},
Chapter II, Proposition 1.8), it follows that for every $x\geq0$,
%
\begin{equation}
\label{eqpropW} \P(Z>x)\leq2\P\biggl[\sup_{t\in[0,1]}W(t)>\frac{x}{2}
\biggr] \leq2\exp\biggl(-\frac{x^2}{8} \biggr).
\end{equation}
In particular, for every $\gamma\in(0,1)$, applying (\ref{eqpropW})
with $x=2\sqrt{2\log(2/\gamma)}$ implies that $q(\gamma
)\leq
2\sqrt{2\log(2/\gamma)}$ and completes the proof.
\end{pf*}


\subsection{\texorpdfstring{Proof of Theorem \protect\ref{theopowerslope}}
{Proof of Theorem 3.1}}
We first prove the following lemma.
%
\begin{lemm}\label{lempower}
Assume ${\cal C}_n$ is any finite collection of subintervals of $[0,1]$ and
%
\begin{equation}
\label{eqcspower} \max_{I\in{\cal C}_n}\sqrt{\frac{n}{|I|}}
\sup_{t\in I} \bigl(\hat F^I(t)-F(t) \bigr)\geq2\sqrt2 \biggl(
\sqrt{\log\biggl(\frac{2|{\cal
C}_{n}|}{\alpha} \biggr)}+\sqrt{\log\biggl(\frac{2}{\beta}
\biggr)} \biggr)
\end{equation}
for some $\alpha$ and $\beta$ in $(0,1)$. Then, (\ref{eqpower}) holds.
\end{lemm}
\begin{pf}
Let $I$ $\in{\cal C}_n$ achieving the maximum in (\ref{eqcspower}) and
$\epsilon=\sqrt{{4\log(2/\beta)}/{\log(2|{\cal C}_{n}|/\alpha)}}$.
It follows from the definition of $S_n^I$, Lemma \ref{lemmajq}
and (\ref{eqcspower}) that
\[
\P_f \Bigl[\max_{I\in{\cal C}_n}S_n^I>s_{\alpha,n}
\Bigr] \geq\P_f \biggl[\sqrt{\frac{n}{|I|}}\sup_{t\in I}
\bigl(\hat F^I_n(t)-F_n(t) \bigr)>2\sqrt{2
\log\biggl(\frac{2|{\cal
C}_{n}|}{\alpha
} \biggr)} \biggr].
\]
Since $\hat F^I_n\geq F_n$ and $\hat F^I\geq F$ on $I$, the triangle
inequality yields
\[
\sup_{t\in I} \bigl(\hat F^I_n(t)-F_n(t)
\bigr)\geq\sup_{t\in I} \bigl(\hat F^I(t)-F(t) \bigr) -
\sup_{t\in I} \biggl|\hat F^I_n(t)-\hat
F^I(t)-\frac{1}{\sqrt
n}W(t)\biggr|.
\]
Besides, with $I=[a,b]$, we have for every $t\in I$:
\[
\biggl|\hat F^I_n(t)-\hat F^I(t)-
\frac{1}{\sqrt n}W(t) \biggr| \leq\biggl|\hat F^I_n(t)-\hat
F^I(t)-\frac{1}{\sqrt n}W(a) \biggr|+\frac
{1}{\sqrt n}\bigl|W(t)-W(a)\bigr|.
\]
But $\hat F^I_n(t)-W(a)/{\sqrt n}$ is the least concave majorant at
time $t$ of the process
$ \{F_n(u)-W(a)/{\sqrt n} \}_{u\in I}$,
so it follows from Lemma 2.2 in \cite{dt03} that
\begin{eqnarray*}
\sup_{t\in I} \biggl|\hat F^I_n(t)-\hat
F^I(t)-\frac
{1}{\sqrt
n}W(a) \biggr| &\leq&\sup_{t\in I}
\biggl|F_n(t)-F(t)-\frac{1}{\sqrt n}W(a) \biggr|
\\
&\leq&\frac{1}{\sqrt n}\sup_{t\in I}\bigl|W(t)-W(a)\bigr|.
\end{eqnarray*}
Hence,
\[
\sup_{t\in I} \biggl|\hat F^I_n(t)-\hat
F^I(t)-\frac{1}{\sqrt
n}W(t) \biggr|\leq\frac{2}{\sqrt n}
\sup_{t\in I}\bigl|W(t)-W(a)\bigr|.
\]
By scaling, the right-hand side is distributed like $2\sqrt{|I|/n}\sup
_{t\in[0,1]}|W(t)|$,
so combining all previous inequalities leads to
\[
\P_f \Bigl[\max_{I\in{\cal C}_n}S_n^I>s_{\alpha,n}
\Bigr]\geq\P_f \biggl[2\sup_{t\in[0,1]}\bigl|W(t)\bigr|<\epsilon\sqrt{2\log
\biggl(\frac{2|{\cal
C}_{n}|}{\alpha} \biggr)} \biggr].
\]
We conclude with the same arguments as in (\ref{eqpropW}) that
\[
\P_f \Bigl[\max_{I\in{\cal C}_n}S_n^I>s_{\alpha,n}
\Bigr]\geq1-2\exp\biggl(-\frac{\epsilon^2}{4}\log\biggl(\frac{2|{\cal
C}_{n}|}{\alpha
}
\biggr) \biggr).
\]
By definition of $\epsilon$, the right-hand side is $1-\beta$, hence
inequality (\ref{eqpower}).
\end{pf}

Let us turn now to the proof of Theorem \ref{theopowerslope} and
recall that we restrict ourselves without loss of generality to the
case $\sigma=1$. Assume (\ref{eqpowerslope}) for some $x,y\in[0,1]$
such that $y-x\geq2/n$, and write $I_0=[x,y]$. The linear function
\[
t\mapsto F(x)+(t-x)\frac{F(y)-F(x)}{y-x}=F(t)+(t-x) ({\bar f}_{xy}-{\bar
f}_{xt}),
\]
where ${\bar f}_{xy}$ is defined by (\ref{eqbarf}), coincides with\vspace*{2pt} $F$
at the boundaries $x$ and $y$ of the interval $I_0$. So this function
is below $\hat F^{I_0}$ on $I_0$ and we obtain
\[
\frac{1}{\sqrt{|I_0|}}\sup_{t\in I_0} \bigl(\hat
F^{I_0}(t)-F(t) \bigr) \geq\sup_{t\in[x,y]}\frac{t-x}{\sqrt{y-x}} (\bar
f_{xy}-\bar f_{xt} ) \geq C(\alpha,\beta)\sqrt
\frac{\log n}{n}.
\]
Now, let $I$ be the smallest interval in ${\cal C}_n$ containing $I_0$.
Since $|I_0|\geq2/n$, we have $|I|\leq|I_0|+2/n\leq2|I_0|$.
Moreover, $\hat F^I\geq\hat F^{I_0}$ on $I_0$, so
%
\begin{equation}
\label{eqmino2} \frac{1}{\sqrt{|I|}}\sup_{t\in I} \bigl(\hat
F^{I}(t)-F(t) \bigr) \geq\frac{1}{\sqrt{2|I_0|}}\sup_{t\in I_0}
\bigl(\hat F^{I_0}(t)-F(t) \bigr) \geq C(\alpha,\beta)\sqrt
\frac{\log n}{2n}.
\end{equation}
Since $|{\cal C}_{n}|=n(n+1)/2$, it follows that (\ref{eqcspower})
holds provided that $C(\alpha,\beta)$ is large enough, so Theorem
\ref
{theopowerslope} follows from Lemma \ref{lempower}.


\subsection{A useful lemma}

In the case $f$ is not non-decreasing such that $f$ is assumed to be
smooth enough, then Lem\-ma~\ref{lemholder} below may serve as a tool to
prove that condition (\ref{eqpowerslope}) in Theorem \ref
{theopowerslope} is fulfilled.
%
\begin{lemm}\label{lemholder}
Assume $f\notin\cal D$ and
$f(u)-f(v)\leq R(u-v)^s$ for all $u\geq v$,
for some $R>0$ and $s\in(0,1]$. Let $x_0<y_0$ in $[0,1]$ such that
$\rho:=f(y_0)-f(x_0)> 0$. Then, there exist an interval $[x,y]\subset
[x_0,y_0]$ and a real $C(s)>0$ that only depends on $s$ such that
%
\begin{equation}
\label{eqsin0,1]} \sup_{t\in[x,y]}\frac{t-x}{\sqrt{y-x}}(\bar f_{xy}-
\bar f_{xt})\geq C(s)R^{-1/(2s)}\rho^{1+1/(2s)}.
\end{equation}
\end{lemm}
\begin{pf} Let $s\in(0,1]$ and $L\geq1$ be fixed, and let
${\mathcal G}(s,L)$ be the set of integrable functions $g\dvtx [0,1]\to\R$
such that $g(0)=0$, $g(1)=1$, and
\[
g(u)-g(v)\leq L(u-v)^s\qquad\mbox{for all }u\geq v.
\]
Let ${\mathcal G}_0(s,L)$ be the set of functions $g\in{\mathcal
G}(s,L)$ such that $\bar g_{01}\geq1/2$, where for every $x<y$, $\bar
g_{xy}$ is defined as in (\ref{eqbarf}).
We first prove that for all $g\in{\mathcal G}(s,L)$,
%
\begin{equation}
\label{eqholdercg} \sup_{0\leq x<t<y\leq1}\frac{t-x}{\sqrt{y-x}}(\bar
g_{xy}-\bar g_{xt})\geq C(s, L):=\frac{s}{7\times2^{1+s+1/s}}L^{-1/(2s)}.
\end{equation}
%
We first consider the case where $g\in{\mathcal G}_0(s,L)$ and argue by
contradiction. Assume there exists $g\in{\mathcal G}_0(s,L)$ such that
inequality (\ref{eqholdercg}) is not satisfied.
For every integer $k$, let $m_k=\bar g_{02^{-k}}$. Setting $x_k=0$,
$t_k=2^{-k-1}$ and $y_k=2^{-k}$, it follows from our assumption on $g$ that
\[
m_k-m_{k+1}=2^{1+k/2}\frac{t_k-x_k}{\sqrt{y_k-x_k}}(\bar
g_{x_ky_k}-\bar g_{x_kt_k}) < 2^{1+k/2}C(s,L).
\]
But for all integers $k_0\geq0$,
\[
\bar g_{01}=m_0=\sum_{k= 0}^{k_0}(m_k-m_{k+1})+m_{k_0+1}.
\]
By assumption, $\bar g_{01}\geq1/2$ whence
%
\begin{equation}
\label{eq12} \frac12< C(s,L)\sum_{k=0}^{k_0}2^{1+k/2}+m_{k_0+1}
< 7C(s,L)2^{k_0/2}+m_{k_0+1}.
\end{equation}
Since $g(0)=0$, for all integers $k_0\geq0$ we have
\[
m_{k_0+1}=2^{k_0+1}\int_0^{2^{-k_0-1}}
\bigl(g(u)-g(0)\bigr) \,\d u< L2^{-s(k_0+1)}.
\]
From (\ref{eq12}) and the definition of $C(s,L)$, we obtain that for
all integers $k_0\geq0$,
%
\begin{equation}
\label{eq122} \frac12< \frac
{s}{2^{1+s+1/s}}L^{-1/(2s)}2^{k_0/2}+L2^{-s(k_0+1)}.
\end{equation}
In particular, consider
$k_0=\sup\{k \in\mathbb N\dvtx  2^{k/2}\leq
2^{1/s}L^{1/(2s)} \}$,
which is well defined since $s \in(0,1]$ and $L\geq1$. By definition
of $k_0$,
$2^{k_0/2}\leq2^{1/s}L^{1/(2s)}$ and $2^{(k_0+1)/2}>
2^{1/s}L^{1/(2s)}$,
so (\ref{eq122}) implies $s2^{-s}> 1/2$.
This is a contradiction because $s2^{-s}\leq1/2$ for all $s\in(0,1]$.
Hence, (\ref{eqholdercg}) holds for all $s\in(0,1]$, $L\geq1$ and
$g\in{\mathcal G}_0(s,L)$. Now, for every $g\in{\mathcal G}(s,L)$ we
set $\tilde g=g$ if $\bar g_{01}\geq1/2$, and $\tilde g(u):=1-g(1-u)$
otherwise, so that $\tilde g\in{\mathcal G}_0(s,L)$. Noting that
\[
\frac{t-x}{\sqrt{y-x}}(\bar g_{xy}-\bar g_{xt})=
\frac{y-t}{\sqrt{y-x}}(\bar g_{ty}-\bar g_{xy})
\]
for all $x<t<y$, we obtain
\[
\sup_{0\leq x<t<y\leq1}\frac{t-x}{\sqrt{y-x}}(\bar g_{xy}-\bar
g_{xt})= \sup_{0\leq x<t<y\leq1}\frac{t-x}{\sqrt{y-x}}(\bar{\tilde
g}_{xy}-\bar{\tilde g}_{xt})\geq C(s,L),
\]
since (\ref{eqholdercg}) holds for all $s\in(0,1]$, $L\geq1$ and
$g\in{\mathcal G}_0(s,L)$. Hence, (\ref{eqholdercg}) holds for all
$s\in(0,1]$, $L\geq1$ and $g\in{\mathcal G}(s,L)$.

Finally, under the assumptions of Lemma \ref{lemholder}, the function
\[
g(u)=\frac{1}{\rho} \bigl(f\bigl(x_0+(y_0-x_0)u
\bigr)-f(x_0) \bigr),\qquad u\in[0,1],
\]
belongs to ${\mathcal G}(s,L)$ with
$L=R(y_0-x_0)^s/\rho\geq1$. Thus, it follows from (\ref
{eqholdercg}) that there exists $C(s)>0$ only depending on $s$ such that
\begin{eqnarray*}
\sup_{x_0\leq x<t<y\leq y_0}\frac{t-x}{\sqrt{y-x}}(\bar f_{xy}-\bar
f_{xt}) &=&\rho\sqrt{y_0-x_0}
\sup_{0\leq x<t<y\leq1}\frac{t-x}{\sqrt{y-x}}(\bar g_{xy}-\bar
g_{xt})
\\
&\geq& C(s)R^{-1/(2s)}\rho^{1+1/(2s)}.
\end{eqnarray*}
\upqed
\end{pf}


\subsection{\texorpdfstring{Remaining proofs for Section \protect\ref{secpower}}
{Remaining proofs for Section 3}}\label{secremproofs}

\begin{pf*}{Proof of Corollary \ref{corolpowerslope}}
Setting $t=t_n$, $x=t_n-\Delta_n$, $y=t_n+\Delta_n$, we obtain
\begin{eqnarray*}
\frac{t-x}{\sqrt{y-x}}(\bar f_{xy}-\bar f_{xt})&=&
\frac{1}{2\sqrt{2\Delta_n}} \biggl(\int_t^yf(u)\,\d u-\int
_x^tf(u)\,\d u \biggr)
\\
&\geq&\frac{1}{2\sqrt{2\Delta_n}} \biggl(\int_{t-\delta_n}^t
\bigl(f(u+\delta_n)-f(u)\bigr)\,\d u \biggr) \geq\frac{M\delta_n^2\lambda
_n}{2\sqrt{2\Delta_n}},
\end{eqnarray*}
so Corollary \ref{corolpowerslope} follows from Theorem \ref
{theopowerslope}.
\end{pf*}
\begin{pf*}{Proof of Corollary \ref{corolpowerUshaped}}
Since $f$ is convex, we can apply Lemma \ref{lemholder} with $s=1$ and
$R$ defined in Corollary \ref{corolpowerUshaped}.
Therefore, there exist $[x,y]\subset[x_{0},1]$ and $C>0$ such that
%
\begin{equation}
\label{eqUinf} \sup_{t\in[x,y]}\frac{t-x}{\sqrt{y-x}}(\bar f_{xy}-\bar
f_{xt})\geq CR^{-1/2}\rho^{3/2}.
\end{equation}
By change of variable, we have for all $t\in[x,y]$,
%
\begin{equation}
\label{eqdiffslope} \bar f_{xy}-\bar f_{xt} =
\frac{1}{y-x}\int_x^y \biggl(f(v)-f \biggl(
\frac
{v-x}{y-x}(t-x)+x \biggr) \biggr)\,\d v.
\end{equation}
By convexity of $f$ and definition of $R$, this implies that for all
$t\in[x,y]$,
\[
\bar f_{xy}-\bar f_{xt}\leq\frac{R}{y-x}\int
_x^y \biggl(\frac
{v-x}{y-x}(y-t) \biggr)\,\d v
\leq\frac{R}{2}(y-x),
\]
hence
%
\begin{equation}
\label{eqUsup} \sup_{t\in[x,y]}\frac{t-x}{\sqrt{y-x}}(\bar f_{xy}-\bar
f_{xt})\leq\frac R2(y-x)^{3/2}.
\end{equation}
Combining inequalities (\ref{eqUinf}) and (\ref{eqUsup}) proves that
$y-x\geq(2C)^{2/3}\rho/R$.
Therefore, $y-x\geq2/n$ provided $\rho>C_{0}R/n$ for a large enough
$C_{0}$.
From Theorem \ref{theopowerslope}, it follows that (\ref{eqcspower})
holds provided
\[
CR^{-1/2}\rho^{3/2}>C(\alpha,\beta)\sqrt{\frac{\sigma^2\log n}{n}}
\]
for a large enough $C(\alpha,\beta)$, which completes the proof of
Corollary \ref{corolpowerUshaped}.
\end{pf*}
\begin{pf*}{Proof of Theorem \ref{theoallpower}} Inequality (\ref
{eqpowerds}) easily follows from Corollary \ref{corolpowerslope}. Yet,
a detailed proof may be found in \cite{CFNsupp}.

We now turn to the proof of inequality (\ref{eqpowerbhl}). Let $f\in
{\cal F}(s,R)$ for some $s\in(0,2]$ and $R>0$ and define
%
\begin{equation}
\label{eqrhobhl} \rho_n=C'(s,\alpha,
\beta)R^{1/(1+2s)} \biggl(\frac{\log n}{n} \biggr)^{s/(1+2s)},
\end{equation}
where $C'(s,\alpha,\beta)$ is a positive number to be chosen later,
that only depends on $s$, $\alpha$ and $\beta$. Assume $n^s\sqrt{\log
n}\geq R$, $R\geq2^{1+2s}\sqrt{(\log n)/n}$ in case $s>1$, and
$\Delta_1(f)\geq\rho_n$.
The function $f^*$ defined on $[0,1]$ by
\[
f^*(y)=\inf_{x\in[0,y]}f(x)
\]
is non-increasing with $f^*\leq f$, so
\[
\Delta_1(f)\leq\sup_{t\in[0,1]}\bigl(f(t)-f^*(t)\bigr) \leq
\sup_{0\leq x<y\leq1}\bigl(f(y)-f(x)\bigr).
\]
Since $f$ is continuous and $\Delta_1(f)\geq\rho_n$, this shows that
there are $x_0<y_0$ (that may depend on $n$) such that
$f(y_0)=f(x_0)+ \rho_n$.

Consider first the case $s\in(0,1]$. From Lemma \ref{lemholder}, there
exist an interval $[x,y]\subset[x_0,y_0]$ and a positive number $C(s)$
only depending on $s$ such that
%
\begin{equation}
\label{eqsupI0} \sup_{t\in[x,y]}\frac{t-x}{\sqrt{y-x}} (\bar
f_{xy}-\bar
f_{xt} )\geq C(s) \bigl(C'(s,\alpha,\beta)
\bigr)^{1+1/(2s)}\sqrt\frac
{\log n}{n}.
\end{equation}
Since $f\in{\cal F}(s,R)$, formula (\ref{eqdiffslope}) implies
\[
\sup_{t\in[x,y]}\frac{t-x}{\sqrt{y-x}} (\bar f_{xy}-\bar
f_{xt} ) \leq\frac{R}{s+1}(y-x)^{s+1/2}.
\]
%
Combining this with (\ref{eqsupI0}) proves that
%
\begin{equation}
\label{eqy-x} y-x\geq\bigl(C(s) (s+1) \bigr)^{2/(1+2s)}
\bigl(C'(s,\alpha,\beta) \bigr)^{1/s} \biggl(
\frac{\log n}{nR^2} \biggr)^{1/(1+2s)},
\end{equation}
so in particular, $y-x\geq2/n$ provided $n^s\sqrt{\log n}\geq R$ and
$C'(s,\alpha,\beta)$ is sufficiently large. Thanks to (\ref{eqsupI0})
and Theorem \ref{theopowerslope}, we obtain that (\ref{eqpower})
holds in the case $s\in(0,1]$, provided $C'(s,\alpha,\beta)$ is large
enough, hence $\rho(T_n,{\cal F}(s,R),\Delta_1)\leq\rho_n$.

Consider now the case $s\in(1,2]$. Assume that, for a given
$C(s,\alpha,\beta)>0$,
%
\begin{equation}
\label{eqmajf} \sup_{t\in[0,1]}f'(t)\leq C(s,
\alpha,\beta)R^{3/(1+2s)} \biggl(\frac
{\log
n}{n} \biggr)^{(s-1)/(1+2s)},
\end{equation}
since otherwise (\ref{eqpowerds}) immediately allows to conclude.
As $f\in{\cal F}(s,R)$ with some $s>1$, we also have
\[
f(u)-f(v)\leq(u-v)\sup_{t\in[0,1]}f'(t)\qquad\mbox{for all }u\geq v,
\]
where $\sup_tf'(t)>0$. Therefore, it follows from the definition of
$x_0,y_0$ and Lemma \ref{lemholder} that there exist $C>0$ and
$[x,y]\subset[x_0,y_0]$ such that
\[
\sup_{x<t<y}\frac{t-x}{\sqrt{y-x}}(\bar f_{xy}-\bar
f_{xt})\geq C \Bigl(\sup_{t\in[0,1]} f'(t)
\Bigr)^{-1/2}\rho_n^{3/2}.
\]
From (\ref{eqmajf}) and the definition of $\rho_n$, we get
\[
\sup_{t\in[x,y]}\frac{t-x}{\sqrt{y-x}}(\bar f_{xy}-\bar
f_{xt})\geq\frac{C(C'(s,\alpha,\beta))^{3/2}}{\sqrt{C(s,\alpha,\beta
)}}\sqrt\frac
{\log n}{n}.
\]
Similar to (\ref{eqy-x}), and using then (\ref{eqmajf}), one obtains
\begin{eqnarray*}
y-x &\geq&\biggl(\frac{2C}{\sup_tf'(t)\sqrt{C(s,\alpha,\beta)}} \biggr
)^{2/3}C'(s,
\alpha,\beta) \biggl(\frac{\log n}{n} \biggr)^{1/3}
\\
&\geq&(2C )^{2/3}\frac{C'(s,\alpha,\beta)}{C(s,\alpha,\beta
)} \biggl(\frac{\log n}{nR^2}
\biggr)^{1/(1+2s)}.
\end{eqnarray*}
Then, we conclude with the same arguments as in the case $s\in(0,1]$.
\end{pf*}


\subsection{\texorpdfstring{Proofs for Section \protect\ref{secregression}}
{Proofs for Section 4}}\label
{secproofreg}
In the sequel, we denote by $H_{\bar n}$ the function $F_{\bar
n}^{\mathrm{reg}}-\sigma_0G_{\bar n}$, which is continuous and piecewise linear
on $[0,1]$:  $H_{\bar n}$ is the cumulative sum diagram of the points
$f_n(i/{\bar n})$, $1\leq i\leq\bar n$, with equal weights $1/\bar n$.

\begin{pf*}{Proof of Theorem \ref{theolevelreg}}
Assume $f\in\mathcal D$. Then, $H_{\bar n}$ has decreasing slopes, so
it is concave on $[0,1]$, and even linear when $f$ is constant over
$[0,1]$. Since $\sigma_0/\hat\sigma_0=\sigma/\hat\sigma$, we
deduce as
in the proof of Theorem \ref{theolevel} that
%
\begin{equation}
\label{eqleastfavreg} \max_{I\in{\cal C}_{\bar n}}\hat S_{\bar
n}^{\mathrm{reg},I}
\leq(\sigma/\hat\sigma) Z^{\mathrm{reg}}_n
\end{equation}
with equality when $f$ is constant over $[0,1]$.
According to Cochran's theorem, $\bar n \hat\sigma^2/\sigma^2$ is
independent of $Z_n^{\mathrm{reg}}$ and distributed as a non-central chi-square
variable with ${\bar n}$ degrees of freedom and non-centrality parameter
%
\begin{equation}
\label{eqnoncentral} \delta^2_f=\bigl(2
\sigma^2\bigr)^{-1}\sum_{i=1}^{{\bar n}}
\bigl(f(2i/n)-f\bigl((2i-1)/n\bigr) \bigr)^2.
\end{equation}
In particular, when $f$ is constant over $[0,1]$,
\[
\hat\sigma^2= \frac{1}{n}\sum_{i=1}^{n/2}
(\epsilon_{2i}-\epsilon_{2i-1})^2:=\hat
\sigma_\epsilon^2,
\]
where ${\bar n}\hat\sigma_\epsilon^2/\sigma^2$ is a $\chi^2({\bar n})$
variable independent of $Z_n^{\mathrm{reg}}$. Thus, we deduce from (\ref
{eqleastfavreg}) and the stochastic order between non-central
chi-square variables with same degrees of freedom (see, e.g.,
\cite{Johnson}, Chapter 29) that
\begin{eqnarray*}
\P_f \Bigl(\max_{I\in{\cal C}_{\bar n}}\hat S_{\bar n}^{\mathrm{reg},I}
>r_{\alpha,n} \Bigr) 
&\leq&\E_f \bigl[
\P_f \bigl(\hat\sigma^2/\sigma^2 <
\bigl(Z^{\mathrm{reg}}_n /r_{\alpha,n}\bigr)^2 |
Z_n^{\mathrm{reg}} \bigr) \bigr]
\\
&\leq&\E\bigl[\P\bigl(\hat\sigma_\epsilon^2/
\sigma^2 < \bigl(Z^{\mathrm{reg}}_n /r_{\alpha,n}
\bigr)^2 | Z_n^{\mathrm{reg}} \bigr) \bigr]
\\
&\leq&\P\bigl((\sigma/\hat\sigma_\epsilon) Z^{\mathrm{reg}}_n
> r_{\alpha,n} \bigr)
\end{eqnarray*}
with equalities in the case where $f$ is constant over $[0,1]$. Both
$\hat\sigma_\epsilon$ and $Z^{\mathrm{reg}}_n$ are continuously distributed, so
their independence, together with the definition of $r_{\alpha,n}$,
shows that the latter probability is equal to $\alpha$.
\end{pf*}

%
\begin{pf*}{Proof of Theorem \ref{theopowerreg}}
As in Sections \ref{secproofslevel} to \ref{secremproofs}, we may
assume that $\sigma=1$. The line of proof of Theorem \ref
{theopowerreg} is close to that of Theorem \ref{theoallpower}.
Indeed, it relies on the discrete versions of Lemmas \ref{lemmajq}
and \ref{lempower} stated below.
\end{pf*}
%
\begin{lemm}\label{lemmajqreg}
For all $\alpha\in(0,1)$ and $n\geq18\log(2/\alpha)$,
\[
r_{\alpha,n}\leq2\sqrt{6\log\biggl(\frac{2|{\cal C}_{\bar
n}|(2-\alpha
)}{\alpha} \biggr)}.
\]
\end{lemm}
%
\begin{lemm}\label{lempowerreg}
Let $\alpha,\beta$ in $(0,1)$, ${\cal C}_{\bar n}$ be the
collection (\ref{eqCnreg}), $L>0$, $s\in(0,2]$, and $R>0$. Assume that
$n\geq18\log(2/\alpha)$ and either that $f\in{\cal F}(s,R)$ for some
$s\in(0,1]$ and $R\leq n^s$, or that $f\in{\cal F}(s,R,L)$ for some
$s\in(1,2]$ and $L\leq n$. Then, there exists a positive number
$C(\alpha,\beta)$ depending only on $\alpha$ and $\beta$ such that, for
all $f$ satisfying
%
\begin{equation}
\label{eqcspowerreg} \max_{I\in{\cal C}_{\bar n}}\sqrt{\frac{\bar n}{|I|}}
\sup_{t\in
I} \bigl({\hat H_{\bar n}}^I(t)-H_{\bar n}(t)
\bigr) \geq C(\alpha,\beta)\sqrt{\log\bar n},
\end{equation}
it holds that
\[
\P_f \Bigl( \max_{I\in{\cal C}_{\bar n}} \hat
S^{\mathrm{reg},I}_{\bar n} > r_{\alpha,n} \Bigr)\geq1-\beta.
\]
\end{lemm}
Detailed proofs of both lemmas and Theorem \ref{theopowerreg} are
given in \cite{CFNsupp}.

\section*{Acknowledgements}

The third author would like to thank Jean-Christophe L\'eger for a
helpful discussion on the proof of Lemma \ref{lemholder}.

\begin{supplement}
\stitle{Supplement to ``Testing monotonicity via local least concave
majorants''}
\slink[doi]{10.3150/12-BEJ496SUPP} 
\sdatatype{.pdf}
\sfilename{BEJ496\_supp.pdf}
\sdescription{We collect in the supplement \cite{CFNsupp} the most
technical proofs. Specifically, we prove how to reduce to the case
$\sigma=1$, we prove (\ref{eqlowerbound2}) and we provide a detailed
proof for (\ref{eqpowerds}) and all results in Section
\ref{secregression}.}
\end{supplement}


\printhistory

\end{document}